\documentclass[aop,MSNbibl,seceqn,citesort,dvips]{arximspdf}
\usepackage{mathbh}
\usepackage{stfloats}
\usepackage{graphicx}

\doi{10.1214/11-AOP710}
\volume{41}
\issue{4}
\pubyear{2013}
\firstpage{2648}
\lastpage{2681}

\makeatletter

\fnbelowfloat

\newcommand{\underset}[2]{\mathop{#2}_{#1}}
\newcommand{\iint}{\int\!\!\int}
\newcommand{\llbracket}{[\![}
\newcommand{\rrbracket}{]\!]}

\renewcommand{\epsilon}{\varepsilon}
\renewcommand{\o}{\mathrm{o}}
\newcommand{\CC}{{\mathbb C}}
\newcommand{\RR}{{\mathbb R}}
\newcommand{\Prob}{\mathbb{P}}
\newcommand{\E}{\mathbb{E}}
\newcommand{\Var}{\operatorname{Var}}
\newcommand{\Tr}{\operatorname{Tr}}
\newcommand{\U}{\mathrm{U}}
\newcommand{\Gue}{\operatorname{GUE}}
\newcommand{\Cue}{\operatorname{CUE}}
\newcommand{\OO}{\mathrm{O}}
\newcommand{\Ext}{\operatorname{Ext}}
\newcommand{\Po}{\operatorname{Po}}
\newcommand{\Cl}{\operatorname{Cl}}
\newcommand{\spec}{\operatorname{spec}}

\newcommand{\one}{\mathbh{1}}
\newcommand{\dd}{{\mathrm{d}}}
\newcommand{\ii}{{\mathrm{i}}}
\newcommand{\Id}{{\mathrm{Id}}}
\renewcommand{\L}{{\mathrm{L}}}
\newcommand{\convlaw}{\stackrel{\mathrm{law}}{\longrightarrow}}
\newcommand{\convLp}{\stackrel{\L^{p}}{\longrightarrow}}

\renewcommand{\Re}{{\mathfrak{Re}}}

\newtheorem{theorem}{Theorem}[section]
\newtheorem{cor}[theorem]{Corollary}
\newtheorem{lem}[theorem]{Lemma}
\newtheorem{prop}[theorem]{Proposition}

\newproclaim{rem}{Remark}

\makeatother

\begin{document}
\begin{frontmatter}

\title{Extreme gaps between eigenvalues of random~matrices}
\runtitle{Extreme gaps between eigenvalues of random matrices}

\begin{aug}
\author[A]{\fnms{G\'erard} \snm{Ben Arous}\thanksref{t1}\ead[label=e1]{benarous@cims.nyu.edu}}
\and
\author[B]{\fnms{Paul} \snm{Bourgade}\corref{}\ead[label=e2]{bourgade@math.harvard.edu}}
\runauthor{G. Ben Arous and Paul Bourgade}
\affiliation{New York University and Harvard University}
\address[A]{Courant Institute\\
\quad of Mathematical Sciences\\
New York University\\
251 Mercer Street\\
New York, New York 10012\\
USA\\
\printead{e1}} 
\address[B]{Department of Mathematics\\
Harvard University\\
One Oxford Street\\
Cambridge, Massachusetts 02138\\
USA\\
\printead{e2}}
\end{aug}

\thankstext{t1}{Supported by the NSF Grants DMS-08-06180 and
OISE-0730136.}

\received{\smonth{5} \syear{2011}}
\revised{\smonth{8} \syear{2011}}

%
\begin{abstract}
This paper studies the extreme gaps between eigenvalues of random
matrices. We give the joint limiting law of the smallest gaps for
Haar-distributed unitary matrices and matrices from the Gaussian
unitary ensemble. In particular, the $k$th smallest gap, normalized by
a factor $n^{-4/3}$, has a limiting density proportional to
$x^{3k-1}e^{-x^3}$. Concerning the largest gaps, normalized by
$n/\sqrt{\log n}$, they converge in $\L^p$ to a constant for all $p>0$.
These results are compared with the extreme gaps between zeros of the
Riemann zeta function.
\end{abstract}

%

\end{frontmatter}

\section{Introduction}


We address here the following question: what is the asymptotic size and
the limit laws of the smallest and
largest gaps or spacings in the spectra of random matrices?
Typical spacings between eigenvalues of random matrices have been well
understood for invariant ensembles for quite some time. More recently, the
behavior of these typical spacings has even been proved to be universal
for much
larger classes of random matrices~\cite{TV,ESY,BG}.
Much less is known for
atypically large or small spacings. This question was first considered
for the smallest spacings
in the unpublished Ph.D. thesis of Vinson~\cite{Vinson} and raised by Diaconis
in~\cite{Diaconis} for the largest ones. It was also discussed in an
interesting debate during
a conference at the Courant Institute in 2006, in honor of Percy Deift.
We solve here completely
the question of the smallest spacings for the simplest invariant
ensembles, that is, the CUE and the GUE. We give
the scaling and the limit laws for the joint distribution of the
smallest spacings. The answer is simple to
state: it is given by the trivial Poissonian ansatz where the spacings
would be treated as i.i.d. random variables.
The answer we find for the largest spacings is less complete, since we
can only obtain at this point a first-order
approximation which gives the asymptotic size of the largest
spacings\vadjust{\goodbreak}
and not their limit laws. We believe that the
same Poissonian ansatz should work as well for the largest gaps, but
this question is left open. The question of the
universality of the behavior of the extremal spacings is also left open.

Dyson~\cite{Dyson1} showed that the repulsion between eigenvalues of
the Gaussian unitary ensemble (the GUE)
could be described asymptotically in terms of the \textit{determinantal}
point process associated with the sine
kernel
\[
K(x,y)=\frac{1}{\pi}\frac{\sin(x-y)}{x-y}.
\]
For this limiting determinantal point process,
the probability of having no eigenvalue in an interval of size $s$ is
known to be the Fredholm determinant
$
\det(\Id-K_{(0,s)}),
$
where $K_{(0,s)}$ is the convolution operator acting on $\L^2(0,s)$
with kernel~$K$. The density of the spacing between two successive
points is then given by (see~\cite{Mehta})
\[
p_2(s)=\partial_{ss}\det\bigl(\Id-K_{(0,s)}\bigr).
\]
This spacing distribution was shown to appear in a number-theoretic
context: some statistics
produced by Odlyzko~\cite{Odl} presented a correspondence between the
histogram of
normalized gaps between zeros of the Riemann zeta function
and $p_2(s)$.
This gave further evidence for the analogy discovered by Dyson and Montgomery:
they realized that the local dependence of the zeros of $\zeta$,
previously calculated by Montgomery, involved the sine kernel; see
\cite{KatzSarnak1,KeaSna} for an historical account
and other steps of this fruitful analogy.
At the same mean or typical gap scale, a precise analysis of the joint
distribution of the gaps between eigenvalues was performed by
Katz and Sarnak~\cite{KatzSarnak1} and Soshnikov~\cite{Sosh1} for the
Circular Unitary Ensemble (the CUE).

Less attention was paid to eigenvalues statistics at smaller and larger
scales. This paper concerns the extreme gaps.
This study was initiated by Vinson
\cite{Vinson}: he showed that the smallest gap between elements of the
$\Cue$, multiplied by $n^{4/3}$, converges in law
to a random variable with distribution function $e^{-x^3}$, as the size
$n$ of the unitary matrix increases.
In his thesis, similar results for the smallest gap between eigenvalues
of a generalization of the $\Gue$ were obtained.
Vinson also gives interesting heuristics suggesting that the largest
gap between $\Cue$ eigenvalues
should be of order $\sqrt{\log n}/n$, with Poissonian fluctuations
around this limit.
Using a different technique, Soshnikov~\cite{Sosh2} investigated the
smallest gaps for determinantal point processes on the real line, with
a translation invariant kernel:
amongst points included in $[0,L]$, this extreme spacing multiplied by
$L^{1/3}$ converges weakly
to the distribution with distribution function $e^{-x^3}$, as $L\to
\infty$.

Heuristically, the above extreme gaps asymptotics can be obtained using
the known asymptotics~\cite{Dyson1,Mehta} of the spacing distribution
\[
p_2(s)\underset{s\to0}{\sim}\frac{\pi^2}{3} s^2, \qquad
\log p_2(s)\underset
{s\to\infty}{\sim}-\frac{s^2}{8},
\]
and treating the gaps as independent random variables. The difficulty
in obtaining rigorous results
lies in showing that this Poissonian ansatz is asymptotically correct
for the extreme gaps, and in making the above estimates uniform in the
dimension~$n$.

We first consider the \textit{joint} law of the smallest gaps (Theorem
\ref{thmunitary}, Corollary~\ref{corsmallestUnitary})
between eigenvalues of unitary matrices. This relies both on
Soshnikov's method and a convergence of the process of small gaps to a Poisson
point process. The same reasoning equally applies to the small gaps
between eigenvalues from the Gaussian unitary ensemble (Theorem \ref
{thmG}, Corollary~\ref{corsmallestGUE}). The proofs of the small
spacings asymptotics are in Section~\ref{sec2}.

The first-order asymptotics of the largest gaps is then proved.
Concerning unitary random matrices (Theorem~\ref{thmlargestGaps}), this
makes use of two important tools. A~key ingredient, by Deift et al.
\cite{Deift}, is the uniform asymptotics about the probability for a
given arc of the circle to be free of eigenvalues. The proof also
requires the negative correlation property for the event that two
disjoint arcs are free of eigenangles. On account of the $\Gue$
(Theorem~\ref{thmlargestGapsGUE}), we also make a essential use of the
negative correlation, and the large gap probability is evaluated by
comparing the $\Gue$ Fredholm determinant with the unitary one. These
large gaps asymptotics are proved in Section~\ref{sec3}.

The extreme spacings between random eigenvalues are important
quantities for statistical physics, computational mathematics and
number theory. For this reason, Diaconis~\cite{Diaconis} mentions the
open question of maximal spacings, answered in Theorem~\ref{thmlargestGaps}.
After making our results explicit, successively for unitary matrices
and the GUE ensemble, we give applications of our extremal spacings
statistics at the end of this Introduction.

\subsection{The unitary group}
Let $u_n$,
a Haar-distributed (measure $\mu_{\U(n)}$) unitary matrix over $\CC
^n$. Suppose
$u_n$ has eigenvalues $e^{\ii\theta_k}$'s,
with ordered eigenangles
$0<\theta_1<\cdots<\theta_n<2\pi$.
Consider the point process on $\RR^2$,
\[
\chi^{(n)}=\sum_{i=1}^n\delta_{(n^{4/3}(\theta_{i+1}-\theta
_i),\theta_i)}.
\]
Our first result is about the convergence of $\chi^{(n)}$ to a Poisson
point process, thanks to this normalization by $n^{-4/3}$.
%
%
\begin{theorem}\label{thmunitary}
Suppose $u_n\sim\mu_{\U(n)}$.
As $n\to\infty$, the process $\chi^{(n)}$ converges to a Poisson point
process $\chi$ with intensity
\[
\E\chi(A\times I)=\biggl(\frac{1}{24\pi}\int_A u^2\,\dd u\biggr)
\biggl(\int
_I\frac{\dd u}{2\pi}\biggr)
\]
for any bounded Borel sets $A\subset\RR^+$ and $I\subset(0,2\pi)$.
\end{theorem}

The intensity is proportional to $\int_I\dd u$ because of the
rotational invariance of the Haar measure. The corresponding factor
will be less\vadjust{\goodbreak}
trivial in the case of the $\Gue$ ensemble.
Our method to prove Theorem~\ref{thmunitary} relies on the
\textit{s-modified} random point field technique
initiated by Soshnikov~\cite{Sosh1,Sosh2}: one can calculate the
correlation functions of the process obtained by keeping
only the $\theta_k$'s, for which $\theta_{k}+An^{-4/3}$ contains
exactly one other eigenvalue. Contrary to~\cite{Sosh1,Sosh2}, we
do not use the notion of cluster functions, because we characterize the
convergence to Poisson random variables, thanks to
the convergence of the factorial moments; this allows us also to
consider easily nontranslation invariant kernels, like in the $\Gue$ case.
Moreover, Theorem~\ref{thmunitary} gives information about the joint
distribution of the number of gaps taking values in disjoint intervals
(convergence in terms of point processes).
In particular we can compute the limiting joint law of the smallest gaps.

Let $t^{(n)}_1<\cdots<t^{(n)}_k$ be the $k$ smallest eigenangles gaps
[i.e., of the form $|\theta_{i+1}-\theta_{i}|$, where the indexes
are modulo $n$ and $|\theta_{i+1}-\theta_{i}|\in(-\pi,\pi)$]. For the
sake of brevity,
write
\[
\tau^{(n)}_k=(72\pi)^{-1/3}t^{(n)}_k.
\]
The limiting joint law of the $\tau_k$'s is a corollary of Theorem
\ref{thmunitary}.
%
%
\begin{cor}\label{corsmallestUnitary}
For any $0\leq x_1<y_1<\cdots<x_k<y_k$, under the Haar
measure on $\U(n)$,
%
%
\begin{equation}\label{eqnjointLarge}\quad
\Prob\bigl(x_\ell<n^{4/3} \tau^{(n)}_\ell<y_\ell,1\leq\ell\leq k\bigr)
\underset{n\to\infty}{\longrightarrow}(e^{-x_k^3}-
e^{-y_k^3})\prod_{\ell=1}^{k-1}(y_\ell^3-x_\ell^3).
\end{equation}
In particular, the $k$th smallest normalized space $n^{4/3}\tau
^{(n)}_k$ converges in law to~$\tau_k$, with density
\[
\Prob(\tau_k\in\dd x)=\frac{3}{(k-1)!}
x^{3k-1}e^{-x^3}\,\dd x.
\]
\end{cor}

Note that this result, for $k=1$, is proved in Vinson's thesis
\cite{Vinson} by
a different method: he characterizes the number of small gaps as a symmetric
function of the eigenvalues, and computes its moments. It is not clear
how his method can be extended to provide
the joint law of the $k$ smallest gaps.

We now turn to our next question about extreme gaps, that is, the
asymptotic behavior of the largest gaps, which were
guessed by Vinson, based on the supposed asymptotic independence of
distant gaps.
We obtain, thanks to the precise asymptotics of one gap probability,
obtained by the steepest descent method for Riemann Hilbert problems in
\cite{Deift},
and the negative association property of determinantal point processes;
see, for example,~\cite{Liggett}.
Note that both results are posterior to Vinson's thesis.

Consider $\mathcal{T}^{(n)}_1>\mathcal{T}^{(n)}_2>\cdots$ the largest
gaps between successive eigenangles of $u\sim\mu_{\U(n)}$,
that is, of the form $|\theta_{i+1}-\theta_{i}|$, where the indexes
are modulo $n$ and $|\theta_{i+1}-\theta_{i}|\in(-\pi,\pi)$.
Then, as $n$ goes to infinity, the largest gap converges in
$\L^{p}$ to a constant, for any $p>0$,
\[
\frac{n}{\sqrt{32\log n}} \mathcal{T}^{(n)}_1\convLp1.
\]

Actually the above limit holds for all the $\ell_n$ largest gaps if $
\ell_n$ is subpolynomial.
%
%
\begin{theorem}\label{thmlargestGaps}
Let $\ell_n=n^{\o(1)}$ be positive integers. Then for any $p>0$,
\[
\frac{n}{\sqrt{32\log n}} \mathcal{T}^{(n)}_{\ell_n}\convLp1
\]
as $n\to\infty$.\setcounter{footnote}{1}\footnote{A detailed analysis
of the proof gives a speed of convergence: for example, one can show
that for any sequence $a_n=\o(1)$,
\[
(\log n)^{a_n}\biggl(\frac{n}{\sqrt{32\log n}} \mathcal
{T}^{(n)}_{1}-1\biggr)\convLp0.
\]
The problem of the exact fluctuations will be addressed in a future work.}
\end{theorem}

Note that, for independent uniform eigenangles on the unit circle, the
largest gap is of order $(\log n)/n$, more than
in the above theorem, as expected from the repulsion of the eigenvalues
in the determinantal case.

\subsection{The Gaussian unitary ensemble}
Similar results hold for the $\Gue$.
For this ensemble, the distribution of the eigenvalues has density
%
%
\begin{equation}\label{eqnGUE}
\frac{1}{Z_n}e^{-n\sum_{i=1}^n\lambda_i^2/2}\prod_{1\leq i<j\leq
n}|\lambda_i-\lambda_j|^2
\end{equation}
with respect to the Lebesgue product measure, on the simplex $\lambda
_1<\cdots<\lambda_n$.
The empirical spectral distribution $\frac{1}{n}\sum\delta_{\lambda_i}$
converges in probability to the
semicircle law (see, e.g.,~\cite{AGZ})
\[
\rho_{\mathrm{sc}}(x)=\frac{1}{2\pi}\sqrt{(4-x^2)_+}.
\]
Like for the unitary group, we first consider the smallest gaps,
studying the point process
\[
\tilde\chi^{(n)}=\sum_{i=1}^{n-1}\delta_{(n^{4/3}(\lambda
_{i+1}-\lambda
_i),\lambda_i)}\one_{|\lambda_i|<2-\epsilon_0}
\]
for any arbitrarily small fixed $\epsilon_0>0$ (this is a technical
restriction allowing the use of the Plancherel--Rotach asymptotics of
the Hermite polynomials).

%
\begin{theorem}\label{thmG}
As $n\to\infty$, the process $\tilde\chi^{(n)}$ converges to a Poisson
point $\tilde\chi$ process with intensity
\[
\E\tilde\chi(A\times I)=\biggl(\frac{1}{48\pi^2}\int_A u^2\,\dd
u
\biggr)\biggl(\int_I(4-x^2)^2\,\dd x\biggr)
\]
for any bounded Borel sets $A\subset\RR_+$ and $I\subset(-2+\epsilon
_0,2-\epsilon_0)$.\vadjust{\goodbreak}
\end{theorem}

The following corollary about the smallest gaps is an easy consequence
of the previous theorem. As for the unitary group, introduce
$\tilde t^{(n)}_1<\cdots<\tilde t^{(n)}_k$ the $k$ nearest spacings in $I$,
that is, of the form $\lambda_{i+1}-\lambda_{i}$, $1\leq i\leq n-1$,
with $\lambda_i \in I$, $I=[a,b]$, $-2<a<b<2$.
Let
\[
\tilde\tau^{(n)}_k=\biggl(\int_I(4-x^2)^2\,\dd x\Big/(144\pi^2)
\biggr)^{1/3}\tilde t^{(n)}_k.
\]

%
\begin{cor}\label{corsmallestGUE}
For any $0\leq x_1<y_1<\cdots<x_k<y_k$, with the above notations and for
the $\Gue$ ensemble measure (\ref{eqnGUE}),
\[
\Prob\bigl(x_\ell<n^{4/3} \tilde\tau^{(n)}_\ell<y_\ell,1\leq
\ell\leq
k\bigr)
\underset{n\to\infty}{\longrightarrow}(e^{-x_k^3}-
e^{-y_k^3})\prod_{\ell=1}^{k-1}(y_\ell^3-x_\ell
^3).
\]
In particular, the $k$th smallest normalized space $n^{4/3}\tilde\tau
^{(n)}_k$ converges in law to~$\tau_k$, with density
\[
\Prob(\tau_k\in\dd x)=\frac{3}{(k-1)!}
x^{3k-1}e^{-x^3}\,\dd x.
\]
\end{cor}
%
%
\begin{cor}\label{corsmallestlocation}
Let $\inf$ be the index of the smallest gap between eigenvalues of the
$\Gue$ in a compact subset $I\subset(-2,2)$ with nonempty interior
\[
\lambda^{(n)}_{\inf+1}-\lambda^{(n)}_{\inf}=\inf\bigl\{\lambda
^{(n)}_{i+1}-\lambda^{(n)}_{i}|\lambda^{(n)}_{i}\in I\bigr\}.
\]
As $n\to\infty$, $\lambda^{(n)}_{\inf}$ converges weakly to the
probability measure with density proportional to
\[
(4-x^2)^2\one_{x\in I}.
\]
\end{cor}

We now turn to the largest gaps for the $\Gue$ ensemble. The result is
completely different inside the bulk and on the edge. Indeed, for eigenvalues
strictly inside the support of the limiting measure, the maximal
spacings have order
$\sqrt{\log n}/n$ (see the following Theorem~\ref{thmlargestGapsGUE}),
while the eigenvalues on the border have an average distance of higher
order, $n^{-2/3}$: for any $k$,
\[
n^{2/3}(\lambda_{n}-2,\ldots,\lambda_{n-k}-2)
\]
converges weakly as $n\to\infty$ to a multivariate Tracy--Widom
distribution; see, for example,~\cite{AGZ}.
Strictly inside the bulk, the result is analogous to the circular case,\vadjust{\goodbreak}
the only difference being the normalization, due to the average density
of eigenvalues.
Let ${\mathcal{\tilde T}}^{(n)}_1>{\mathcal{\tilde T}}^{(n)}_2>\cdots$
be the largest gaps of type $\lambda_{i+1}-\lambda_i$ with $\lambda
_i\in I$, a compact subset of $(-2,2)$ with nonempty interior.
%
%
\begin{theorem}\label{thmlargestGapsGUE}
Let $\ell_n=n^{\o(1)}$ be positive integers. Then for any $p>0$,
\[
\Bigl(\inf_{x\in I}\sqrt{4-x^2}\Bigr) \frac{n}{\sqrt{32\log
n}}
\mathcal{\tilde T}^{(n)}_{\ell_n}\convLp1.
\]
\end{theorem}

\subsection{\texorpdfstring{The $\zeta$ zeros}{The zeta zeros}}
When seen in a window of size proportional to the average gap, the
spacings between the zeros of Dirichlet L-functions are distributed
like particles of a determinantal point process with sine kernel; this
is the Montgomery--Odlyzko law~\cite{Montg}. Here we want to discuss
the accuracy of this analogy when looking at rare events, the extreme
gaps between the zeta zeros, relying on Theorems~\ref{thmunitary} and
\ref{thmlargestGaps}.

Due to the availability of many numerical data, we focus on the Riemann
zeta function
$\zeta(s)=\sum_{n=1}^\infty1/n^s$ ($\Re(s)>1$), which admits an
analytic continuation to $\CC-\{1\}$.
Let $1/2\pm\ii t_k$ be the nontrivial zeta zeros, $\gamma_i=\Re(t_i)$,
with $0<\gamma_1\leq\allowbreak\gamma_2\leq\cdots.$
Then
\[
\tilde\gamma_i=\frac{\gamma_{i+1}-\gamma_i}{2\pi}\log
\biggl(\frac{\gamma
_i}{2\pi}\biggr)
\]
has an average value 1.
The quantity
\[
\lambda=\limsup_{i\to\infty}\tilde\gamma_i
\]
has been widely studied.
Conditionally to the generalized Riemann hypothesis, the best known
result is $\lambda>3.0155$~\cite{Bui}.
From the $\Gue$ hypothesis for the zeta zeros, it is expected that
$\lambda=\infty$. However, to the best of our knowledge,
more precise conjectures about the growth speed of large gaps between
zeta zeros were not proposed.
From Theorem~\ref{thmlargestGaps}, amongst $n$ successive gaps with
fermionic repulsion, the maximal gap has size about
$\sqrt{32\log n}/(2\pi)$ times the average gap, suggesting
\[
\sup_{m\leq k\leq m+n}\tilde\gamma_i\underset{n\to\infty}{\sim
}\frac
{\sqrt{32\log n}}{2\pi},
\]
in particular
\[
\limsup_{i\to\infty}\sqrt{\frac{\log\gamma_i}{32}}(\gamma
_{i+1}-\gamma_i)=1.
\]
Odlyzko's numerical data~\cite{Odl} give 3.303 for the maximal value of
$\tilde\gamma_i$, $1\leq i\leq n=10^6$, while $\frac{\sqrt{32\log
n}}{2\pi}=3.346$, giving a difference of $1\%$ with the observed gaps.
%
%
\begin{figure}

\includegraphics{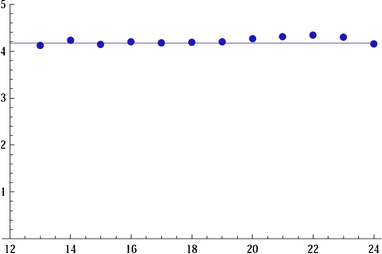}

\caption{Large gaps between the zeta zeros.}
\label{figmaximal}
\end{figure}
Further tests can be performed at distinct heights along the critical
axis, thanks to numerical data of Gourdon~\cite{Gourdon}: he computed
$n=2\times10^9$ successive zeta zeros at height $10^k$ along the
critical axis for each $k\in\llbracket13,24\rrbracket$.
The extreme normalized gaps are given on the joint Figure~\ref{figmaximal}, where the
expectation from our random matrices result is the straight line $\frac
{\sqrt{32\log n}}{2\pi}=4.166$.
For example, amongst the $2\times10^9$ gaps following the height
$10^{24}$, $\sup\tilde\gamma_i=4.158$, that is, a difference of
$0.2\%$
with the expected value.

Concerning the smallest gaps, does the Poisson intensity $\frac
{1}{24\pi
}u^2\,\dd u$ from Theorem~\ref{thmunitary} appear in the context of the
zeta zeros?
Note $\tilde\theta_i=n(\theta_{i+1}-\theta_i)/(2\pi)$ the normalized
gaps. We know that, as $n\to\infty$, the set of gaps $\{2\pi n^{1/3}
\tilde\theta_i,1\leq i\leq n\}$ converges weakly
to a Poisson point process with intensity $\frac{1}{24\pi}u^2\,\dd u$.
We therefore expect that, as $n\to\infty$,
\[
2\pi n^{1/3}\{\tilde\gamma_i,1\leq i\leq n\}
\]
converges to a Poisson point process with the same intensity.
%
%
\begin{figure}[b]

\includegraphics{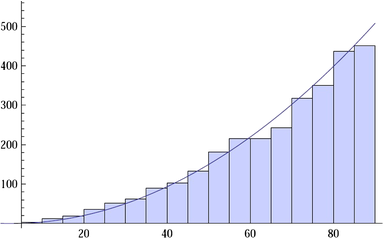}

\caption{Small gaps between the zeta zeros.}\label{fig2}
\end{figure}
The joint Figure~\ref{fig2} gives the histogram of the $3000$ smallest gaps,
normalized as previously, amongst the $n=10^{13}$ first zeta zeros,
based again on numerical data from~\cite{Gourdon}. More precisely, the
histogram gives the 3000 smallest values of
$2\pi10^{13/3}\{\tilde\gamma_i,1\leq i\leq10^{13}\}$.
The straight line is the function $5\times\frac{1}{24\pi}u^2$ (the step
of the histogram is 5). This presents a good relevance of the $\Gue$
hypothesis for the Riemann zeta function, even at the scale of rare
events, here
the extreme spacings.

\subsection{Diagonalization speed with the Toda flow}
The most classical method to diagonalize a matrix is the well-known QR
algorithm. In the case of Hermitian matrices, an alternative approach
was proposed by
Deift et al.~\cite{DNT},\footnote{Their approach is enounced for
symmetric matrices, and naturally extends to the Hermitian case.} based
on the isospectral property
of the Toda flow. More precisely, given a $n\times n$ Hermitian matrix
$M$, the first step is to reduce it in a tridiagonal form $T$ (which is
a robust and fast operation), conjugating with successive Householder
reflections,
\[
T=
\left[
\matrix{
a_1&b_1&&\cr
b_1&a_2&\ddots&\cr
&\ddots&\ddots&b_{n-1}\cr
&&b_{n-1}&a_n}
\right],
\]
keeping the same spectrum as $M$. After such a tridiagonalization, the
$a_i$ and $b_i$'s are real.
The above matrix is important in the analysis of the Toda lattice:
Flaschka, H\'enon and Manakov proved independently in the 1970s that
the following evolution of $n$ particles on a line ($x_0=-\infty$,
$x_{n+1}=+\infty$, $1\leq k\leq n$),
\[
\ddot{x}_k=e^{x_{k-1}-x_{k}}-e^{x_{k}-x_{k+1}}
\]
is an integrable system. More precisely, after the change of variables
\[
\cases{
a_i=-\dot{x}_i/2,\vspace*{2pt}\cr
b_i=\frac{1}{2}e^{(x_i-x_{i+1})/2},}
\]
the differential equation takes the Lax pair form
\[
\frac{\dd T}{\dd t}=S T-T S,
\]
where
\[
S=
\left[
\matrix{
0&b_1&&\cr
-b_1&0&\ddots&\cr
&\ddots&\ddots&b_{n-1}\cr
&&-b_{n-1}&0}
\right].
\]
In particular and importantly, the spectrum of $T(t)$ does not depend
on time.
Moser proved that $\dot{x}_k(t)=\lambda_k+\o(1)$,
$x_k(t)=\lambda_k t+\mu_k+\o(1)$ as $t\to\infty$, with \mbox{$\lambda
_1<\cdots<\lambda_n$}.
This implies that $b_k(t)$ converges to $0$, hence $T$ converges to a
diagonal matrix, whose entries give the eigenvalues of $M=T(0)$.
Deift~\cite{DeiftOpenPbs} asked about the speed of convergence of the
Toda flow
till its equilibrium. More precisely, for a given $\epsilon>0$, what is
the necessary time $t$ such that
$b_k(t)<\epsilon$
for all $1\leq k\leq n-1$? As
\[
b_k(t)=e^{({1}/{2})(\lambda_k-\lambda_{k+1})t+({1}/{2})(\mu
_k-\mu
_{k+1})+\o(1)},
\]
the speed of convergence to the spectrum is governed by the minimal gap
between eigenvalues. A good choice for a \textit{typical} Hermitian matrix
is a matrix from the $\Gue$, with independent (up to symmetry)
complex Gaussian entries (of variance~2 on the diagonal, 1 elsewhere).
From Corollary~\ref{corsmallestGUE}, the minimal gap for such matrices
is of order
\[
\sqrt{n} n^{-4/3}=n^{-5/6}.
\]
For a given required precision $\epsilon$, the
Toda flow time necessary to evaluate the eigenvalues is expected to
grow as $n^{5/6}$ with the dimension.\footnote{Note that more precise estimates
would need the asymptotics of $\mu_k-\mu_{k+1}$ as a function of $n$.}


\section{Small gaps}\label{sec2}

\subsection{Convergence to Poisson point processes}
We first review some results about convergence of processes, here with
values in $\RR^2$ and a finite number of atoms.
A point process $\chi^{(n)}=\sum_{i=1}^{k_n}\delta_{X_{i,n}}$ is said
to converge in distribution
to $\chi=\sum_{i}\delta_{X_{i}}$ is, for any bounded continuous
function $f$,
\[
\chi^{(n)}(f)=\sum_{i=1}^{k_n}f(X_{i,n})\convlaw\sum_{i}f(X_i).
\]
%

To show the convergence of $\chi^{(n)}$ to a Poisson point process
$\chi$,
we only need to show the convergence in law of $\chi^{(n)}(A,I)$ to
$\chi(A,I)$ for all
bounded intervals $A$ and $I$, the independence
for disjoint $A\times I$'s being
an automatic consequence: this is a very practical property of Poisson
point processes,
detailed in Proposition~\ref{propfacMom} below.
Moreover, the convergence of $\chi^{(n)}(A,I)$ will be shown, thanks to
the convergence of the factorial moments
to those of a Poisson random variable.
This is particularly adapted to our situation because, as we will see,
the correlation functions of point processes are defined through
factorial moments and are explicit in the case of determinantal point processes.

Note that this is the same technique employed in~\cite{BK}, where the
following result
is given, in the case of point processes with values in $\RR$.
%
%
\begin{prop}\label{propfacMom}
Let $\chi^{(n)}=\sum_{i=1}^{k_n}\delta_{X_{i,n}}$ be a sequence of
point processes on~$\RR^2$, and $\chi$ a Poisson point process on\vadjust{\goodbreak}
$\RR^2$ with intensity $\mu$ having no atoms (and $\sigma$-finite).
Assume that for any bounded intervals $A$ and $I$ and all positive
integers $k\geq1$,
%
%
\begin{equation}\label{eqnFacMomPoisson}
\lim_{n\to\infty}\E\biggl(\frac{\chi^{(n)}(A\times I)!}{(\chi
^{(n)}(A\times I)-k)!}\biggr)=\mu(A\times I)^k.
\end{equation}
Then the sequence of point processes $\chi^{(n)}$ converges in
distribution to $\chi$.
\end{prop}
\begin{pf}
We need to check the three conditions of the following Theorem~\ref
{thmKallenberg} by Kallenberg, which is written here in the specific
case $\mathfrak{S}=\RR^2$, $\mathcal{U}$ and $\mathcal{J}$ the set of
compact rectangles $A\times I$ ($A$ and $I$ intervals),
with the notation from~\cite{Kallenberg}.
Conditions (1), (2) and (3) will be verified if
%
%
\begin{equation}\label{eqnconvlaw}
\chi^{(n)}(A\times I)\convlaw\chi(A\times I).
\end{equation}
Indeed, if (\ref{eqnconvlaw}) holds, for any $t$,
\[
\lim_{n\to\infty}\Prob\bigl(\chi^{(n)}(A\times I)>t\bigr)=
\Prob\bigl(\chi(A\times I)>t\bigr),
\]
and as $\chi(A\times I)$ is almost surely finite, this goes to $0$ as
$t\to\infty$
(is is a consequence of the dominated convergence theorem); this proves
part (1). Equations (2) and (3) are also
consequences of the above convergence in law.

To prove (\ref{eqnconvlaw}), as $\chi(A\times I)$ is a Poisson
random variable,
the convergence of all the moments is sufficient. A moment is a finite
linear combination of factorial moments, which concludes the proof.
\end{pf}
%
%
\begin{theorem}[(Kallenberg~\cite{Kallenberg})]\label{thmKallenberg}
Let $\chi$ be a point process on $\RR^2$, and assume $\chi$ is almost
surely simple
(i.e., the atoms of the measure $\chi$ all have weight 1 almost surely).
Then $\chi^{(n)}$ converges weakly to $\chi$ if and only if the three
following conditions are satisfied
for any compact intervals $A$ and $I$ in~$\RR$:
\begin{longlist}[(2)]
\item[(1)] $\lim_{t\to\infty}\lim_{n\to\infty}\Prob(\chi
^{(n)}(A\times
I)>t)=0$;
\item[(2)] $\lim_{n\to\infty}\Prob(\chi^{(n)}(A\times I)=0
)=\Prob
(\chi(A\times I)=0)$;
\item[(3)] $\limsup_{n\to\infty}\Prob(\chi^{(n)}(A\times
I)>1)\leq
\Prob(\chi(A\times I)>1)$.
\end{longlist}
\end{theorem}

\subsection{Correlation functions}
References for the main properties of correlation functions of
determinantal point processes are~\cite{JohSurvey} and~\cite{SoshSurvey}.
We follow this last survey to present the notions used in the following.
If $\chi=\sum_i\delta_{X_i}$ is a simple point process on a complete
separate metric space $\Lambda$, consider the point process
%
%
\begin{equation}\label{eqnmultdimPcess}
\Xi^{(k)}=\sum_{X_{i_1},\ldots,X_{i_k}\ \mathrm{all}\
\mathrm{distinct}}\delta_{(X_{i_1},\ldots,X_{i_k})}
\end{equation}
on $\Lambda^k$. One can define this way a measure $M_k$ on $\Lambda
^k$ by
\[
M^{(k)}(\mathcal{A})=\E\bigl(\Xi^{(k)}(\mathcal{A})\bigr)\vadjust{\goodbreak}
\]
for any Borel set $\mathcal{A}$ in $\Lambda^k$. Most of the time, there
is a natural measure $\lambda$ on~$\Lambda$, in our cases $\Lambda=\RR$ or $(0,2\pi)$, and $\lambda$ is
the Lebesgue measure.
If $M^{(k)}$ is absolutely continuous with respect to $\lambda^k$,
there exists a function $\rho_k$ on $\Lambda^k$ such that for any Borel
sets $B_1,\ldots,B_k$ in $\Lambda$
\[
M^{(k)}(B_1\times\cdots\times B_k)=\int_{B_1\times\cdots\times
B_k}\rho
_k(x_1,\ldots,x_k)\,\dd\lambda(x_1)\cdots\dd\lambda(x_k).
\]
Hence one can think about $\rho_k(x_1,\ldots,x_k)$ as the asymptotic
(normalized) probability of having exactly one particle in
neighborhoods of the $x_k$'s. More precisely, under suitable smoothness
assumptions, and for distinct points
$x_1,\ldots,x_k$ in $\Lambda=\RR$,
\[
\rho_k(x_1,\ldots,x_k)=\lim_{\epsilon\to0}\frac{1}{\prod
_{j=1}^k\lambda
(x_j,x_j+\epsilon)}\Prob\bigl(\chi(x_i,x_i+\epsilon)=1,1\leq i\leq
k\bigr).
\]
Note that $\rho_k$ is not a probability density.
Moreover, specifically, if $B_1,\ldots, B_\ell$ are disjoint in
$\Lambda
$ and $n_1+\cdots+n_\ell=k$,
\[
M^{(k)}(B_1^{n_1}\times\cdots\times B_\ell^{n_\ell})=\E\Biggl(\prod
_{i=1}^\ell\frac{(\chi(B_i))!}{(\chi(B_i)-n_i)!}\Biggr);
\]
see~\cite{JohSurvey} for a proof. In particular,
%
%
\begin{equation}\label{eqnFacMomCorrelation}\qquad
M^{(k)}(B^k)=\E\biggl(\frac{(\chi(B))!}{(\chi(B)-k)!}\biggr)=\int
_{B^k}\rho_k(x_1,\ldots,x_k)\,\dd\lambda(x_1)\cdots\dd\lambda(x_k).
\end{equation}
Note the analogy with formula (\ref{eqnFacMomPoisson}) we want to prove.
For unitary matrices or the $\Gue$ ensemble, our method to prove
convergence of small spacings counting measures is the same:
\begin{itemize}
\item For given compact intervals $A$ and $I$, consider the modified
process obtained from $\xi^{(n)}=\sum\delta_{X_{i,n}}$
by keeping only the points $X_{i,n}$ in $I$ such that $\chi
^{(n)}(X_{i,n}+A n^{-4/3})=1$.
\item Show that the correlation function $\tilde\rho
^{(n)}_k(x_1,\ldots,x_n)$ of this new process
uniformly converges to $\mu(A\times I)^k$. This is possible, thanks to
the determinantal aspect of
$\xi^{(n)}$ and the
Hadamard--Fischer inequality, Lemma~\ref{lemFischer}.
\item Conclude that the factorial moments converge to those of the
expected Poisson random variables,
thanks to (\ref{eqnFacMomPoisson}) and (\ref{eqnFacMomCorrelation}).
\end{itemize}

For the smallest gaps asymptotics, the following inequality will be
repeatedly used.
A concise proof can be found in~\cite{HornJohnson}.
%
%
\begin{lem}\label{lemFischer}
Let $M$ be a positive-definite $n\times n$ (Hermitian) matrix. For any
$\omega\subset\llbracket1,n\rrbracket$, let
$M_\omega$ (resp., $M_{\overline{\omega}}$) be the submatrix of $M$
using rows and columns numbered in $\omega$
(resp., $\llbracket1,n\rrbracket\slash\omega$). Then
\[
\det(M)\leq\det(M_\omega)\det(M_{\overline{\omega}}).
\]
\end{lem}

\subsection{The unitary group}
We begin with the proof of Theorem~\ref{thmunitary}.
We know that for a unitary matrix $u_n\sim\mu_{\U(n)}$, the density of
the eigenangles
$0\leq\theta_1<\cdots<\theta_n<2\pi$, with respect to the Lebesgue
measure on the corresponding simplex is
\[
\frac{1}{(2\pi)^n}\prod_{j<k}|e^{\ii\theta_j}-e^{\ii\theta_k}|^2.
\]
Moreover, a remarkable fact about the point process $\sum\delta
_{e^{\ii
\theta_k}}$ is that it is determinantal: all its
correlation functions $\rho^{\U(n)}_k$, $1\leq k\leq n$, are
determinants based on the same kernel,
\[
\rho^{\U(n)}_k(\theta_1,\ldots,\theta_n)=\det_{k\times k}
\bigl(K^{\U
(n)}(\theta_i-\theta_j)\bigr),\qquad K^{\U(n)}(\theta)=\frac{1}{2\pi
}\frac
{\sin(n\theta/2)}{\sin(\theta/2)}.
\]
This classical property relies on Gaudin's lemma; see~\cite{Mehta}.
In the following, for any bounded interval $A\subset\RR^+$, we write
$A_n=n^{-4/3}A$.
We want to show that for an interval $I\subset(0,2\pi)$,
\[
\chi^{(n)}(A\times I)\convlaw\Po(\lambda)
\]
with $\lambda=(\frac{1}{48\pi^2}\int_A u^2\,\dd u)
(\int
_I\dd u)$.
We consider the point process
\[
\xi^{(n)}=\sum_{i=1}^n\delta_{\theta_i}
\]
and its thinning $\tilde\xi^{(n)}$
obtained from $\xi^{(n)}$ by only keeping the eigenangles $\theta_k$
for which $\xi^{(n)}(\theta_k+A_n)=1$.
The following lemma means
that $\chi^{(n)}(A\times I)$ is properly estimated by $\tilde\xi
^{(n)}(I)$. It is analogous to Lemma 3 in~\cite{Sosh2}.
%
%
\begin{lem}[(No successive small neighbors)]\label{lemnoSucSmall}
For any interval $I\subset(0,2\pi)$, as $n\to\infty$, $\chi
^{(n)}(A\times I)-\tilde\xi^{(n)}(I)\convlaw0$.
\end{lem}
\begin{pf}
Let $c$ be such that $A\subset(0,c)$, and $c_n=c n^{-4/3}$. If $\one
_{\theta_{i+1}-\theta_i\in A_n}\neq\one_{\xi^{(n)}(\theta_i+A_n)=1}$,
then $\xi^{(n)}(\theta_i+(0,c_n))\geq2$. Hence
\[
\bigl|\chi^{(n)}(A\times I)-\tilde\xi^{(n)}(I)\bigr|\leq\sum_{i=1}^n\one
_{\xi
^{(n)}(\theta_i+(0,c_n))\geq2}\leq\Xi^{(3)}(\mathcal{A}),
\]
where the last inequality comes from the definition (\ref
{eqnmultdimPcess}), where $\mathcal{A}$
is the set of points $(\theta,x_1,x_2)$ with $\theta\in(0,2\pi)$ and
$(x_1,x_2)\in(\theta,\theta+c_n)^2$.
To prove that this positive random variable converges in law to $0$, we
consider its expectation
\[
\int_{0}^{2\pi}\dd\theta
\int_{(\theta,\theta+c_n)^2}\rho^{\U(n)}_3(\theta,x_1,x_2)\,\dd
x_1\,\dd
x_2=2\pi\int_{(0,c_n)^2}
\rho^{\U(n)}_3(0,x_1,x_2)\,\dd x_1\,\dd x_2
\]
and show it goes to 0. Thanks to the multilinearity of the determinant,
\begin{eqnarray*}
\hspace*{-3pt}&&\rho^{\U(n)}_3(0,x_1,x_2)
\\
\hspace*{-3pt}&&\quad=
{\fontsize{10.7pt}{11pt}\selectfont{\left|
\matrix{
K^{\U(n)}(0)&K^{\U(n)}(x_1)-K^{\U(n)}(0)&K^{\U(n)}(x_2)-K^{\U
(n)}(0)\vspace*{2pt}\cr
K^{\U(n)}(x_1)&K^{\U(n)}(0)-K^{\U(n)}(x_1)&K^{\U(n)}(x_1-x_2)-K^{\U
(n)}(x_1)\vspace*{2pt}\cr
K^{\U(n)}(x_2)&K^{\U(n)}(x_1-x_2)-K^{\U(n)}(x_2)&K^{\U(n)}(0)-K^{\U
(n)}(x_2)}
\right|}}.
\end{eqnarray*}
As $|K^{\U(n)}|_\infty=O(n)$ and $|{K'^{\U(n)}}|_\infty=O(n^2)$, the
first column of this determinant is $O(n)$, and the two others are
$O(n^2 c_n)$.
Thus
\[
\rho^{\U(n)}_3(0,x_1,x_2)=O(n^{7/3}).
\]
The integration domain is $c_n^2=O(n^{-8/3})$, concluding the proof.
\end{pf}

Let $\tilde\rho^{\U(n)}_k(\theta_1,\ldots,\theta_k)$, $k\geq0$,
be the
correlation functions of the point process~$\tilde\xi^{(n)}$.
If, for any $k\geq1$, the convergence of the factorial moment
%
%
\begin{eqnarray}\label{eqnconvFacMom}
\E\biggl(\frac{(\tilde\xi^{(n)}(I))!}{(\tilde\xi
^{(n)}(I)-k)!}\biggr)&=&
\int_{I^k}\tilde\rho^{\U(n)}_k(\theta_1,\ldots,\theta_k)\,\dd
\theta_1\cdots
\dd\theta_k\nonumber\\[-8pt]\\[-8pt]
&\displaystyle \underset{n\to\infty}{\longrightarrow}&\biggl(\frac{1}{24\pi}\int_A
u^2\,\dd u\biggr)^k\biggl(\int_I\frac{\dd u}{2\pi}\biggr)^k
\nonumber
\end{eqnarray}
can be shown, then Theorem~\ref{thmunitary} will be proved, thanks to
the above Lem\-ma~\ref{lemnoSucSmall}.
The way to show (\ref{eqnconvFacMom}) relies on three steps, to apply
a simple dominated convergence argument:
\begin{itemize}
\item if all $\theta_k$'s are distinct, $\tilde\rho^{\U
(n)}_k(\theta
_1,\ldots,\theta_k)$ converges to $(\frac{1}{48\pi^2}\int_A
u^2\,\dd
u)^k$ as
$n\to\infty$ (Lemma~\ref{lemsimpleConvergence});
\item in the set
%
%
\begin{equation}\label{eqndefOmega}
\Omega^{\U(n)}=\{(\theta_1,\ldots,\theta_k)\in I^k\dvtx\theta_i\notin
\theta
_j+A_n,1\leq i,j\leq k\},
\end{equation}
$\tilde\rho{}^{\U(n)}_k(\theta_1,\ldots,\theta_k)$ is uniformly bounded
(Lemma~\ref{lemuniformBounded});
\item even if $\tilde\rho^{\U(n)}_k(\theta_1,\ldots,\theta_k)$ is not
uniformly bounded in the complement of $\Omega^{\U(n)}$ in $I^k$
($\overline{\Omega}{}^{\U(n)}$),
the contribution to the integral is negligible because the volume of
$\overline{\Omega}{}^{\U(n)}$ decreases sufficiently fast (Lemma \ref
{lemneglectSet}).
\end{itemize}
%
%
\begin{lem}[(Simple convergence)]\label{lemsimpleConvergence}
Let $\theta_1,\ldots,\theta_k$ be distinct elements in~$I^k$. Then
\[
\tilde\rho^{\U(n)}_k(\theta_1,\ldots,\theta_k)\underset{n\to
\infty
}{\longrightarrow}\biggl(\frac{1}{48\pi^2}\int_A u^2\,\dd u\biggr)^k.
\]
\end{lem}
\begin{pf}
First note that, as all the $\theta_k$'s are distinct,
for sufficiently large~$n$, the point\vspace*{1pt} $(\theta_1,\ldots,\theta_k)$ is in
$\Omega^{\U(n)}$; see (\ref{eqndefOmega}). This means
that if $\theta_1,\ldots,\theta_k$ are points of~$\tilde\xi^{(n)}$, the\vadjust{\goodbreak}
point in each of the $\theta_i+A_n$
is not another one of the $\theta_j$'s. This makes the combinatorics easy:
the correlation functions of $\tilde\xi^{(n)}$
can be explicitly given in terms of those of $\xi^{(n)}$, as noted in
\cite{SoshSurvey}, by an inclusion--exclusion argument:
for sufficiently large $n$,
%
%
\begin{eqnarray}\label{eqninclusionExclusion}
&&\tilde\rho^{\U(n)}_k(\theta_1,\ldots,\theta_k)\nonumber\\
&&\qquad=\sum_{m=0}^\infty
\frac
{(-1)^m}{m!}\int_{\theta_1+A_n}\dd x_1\cdots
\int_{\theta_k+A_n}\dd x_k\nonumber\\[-8pt]\\[-8pt]
&&\qquad\quad{}\times\int_{((\theta_1+A_n)\sqcup\cdots\sqcup(\theta_k+A_n))^m}\rho^{\U
(n)}_{2k+m}(\theta_1,x_1,\ldots,\theta_k,x_k,\nonumber\\
&&\qquad\quad\hspace*{173pt}y_1,\ldots,y_m)\,\dd
y_1\cdots
\dd y_m.
\nonumber
\end{eqnarray}
Note that there is no convergence issue here as $\rho^{\U
(n)}_{2k+m}\equiv0$ if $2k+m>n$. We first show that the
term corresponding to $m=0$ in the above sum gives the expected asymptotics.
The determinantal aspect of the process makes things easy.
$\rho^{\U(n)}_{2k}(\theta_1,x_1,\ldots,\theta_k,x_k)$ is a
$2k\times2k$
determinant, and
only the terms in the $2\times2$ diagonal blocks make a significant
contribution (which leads to the idea of asymptotic independence).
More precisely, we write formally
\[
\rho^{\U(n)}_{2k}(\theta_1,x_1,\ldots,\theta_k,x_k)=\det_{1\leq
i,j\leq
k}\pmatrix{
K^{\U(n)}(\theta_i-\theta_j)&K^{\U(n)}(\theta_i-x_j)\cr
K^{\U(n)}(x_i-\theta_j)&K^{\U(n)}(x_i-x_j)}.
\]
As $|x_i-\theta_i|=O(n^{-4/3})$ and $K^{\U(n)}(x)=\sin(n x /2)/\sin
(x/2)$ if $i\neq j$ all terms of the corresponding $2\times2$ above
matrix are $O(1)$. Moreover, the above determinant is unchanged by
subtracting an odd column to the following even column, and then by
subtracting an odd line to the following even line. In this way, the
diagonal $2\times2$ matrices becomes
\begin{eqnarray*}
&&\pmatrix{
K^{\U(n)}(0)&K^{\U(n)}(\theta_i-x_i)-K^{\U(n)}(0)\vspace*{2pt}\cr
K^{\U(n)}(x_i-\theta_i)-K^{\U(n)}(0)&2K^{\U(n)}(0)-K^{\U
(n)}(x_i-\theta
_i)-K^{\U(n)}(\theta_i-x_i)}
\\
&&\qquad=
\pmatrix{
O(n)&O(n^{2/3})\vspace*{2pt}\cr
O(n^{2/3})&O(n^{1/3})},
\end{eqnarray*}
where the last equality\vspace*{1pt} relies on $|K^{\U(n)}|_\infty=O(n)$, $|K'^{\U
(n)}|_\infty=O(n^2)$ and $|K''^{\U(n)}|_\infty=O(n^3)$.
As a consequence, in the expansion of the determinant over all
permutations of $\mathcal{S}_{2k}$,
the terms corresponding to entries only in the diagonal $2\times2$
block matrices have order at most
$n^{({4}/{3})k}$, while all other terms have a strictly lower order
(at most $n^{({4}/{3})k-({2}/{3})}$). As the integration domain of
$\rho^{\U(n)}_{2k}$ is $O(n^{-({4}/{3})k})$, the only permutations
hopefully giving a nonzero limit need to come from the block diagonal
$2\times2$ matrices.
Indeed they give a nontrivial limit: their contribution is exactly
\[
\prod_{i=1}^k\int_{\theta_i+A_n}\rho^{\U(n)}_2(\theta_i,x)\,\dd
x=
\biggl(\frac{1}{(2\pi)^2}\int_{A_n}n^2\biggl(1-\biggl(\frac{\sin(n
x/2)}{n\sin
(x/2)}\biggr)^2\biggr)\,\dd x\biggr)^k.
\]
A simple change of variable $x=n^{4/3}u$ allows us to conclude, thanks
to the easy limit, uniform on compacts,
\[
\frac{1}{(2\pi)^2}n^{2/3}\biggl(1-\biggl(\frac{\sin(n^{-1/3}
u/2)}{n\sin
(n^{-4/3}u/2)}\biggr)^2\biggr)\underset{n\to\infty
}{\longrightarrow}
\frac{1}{48\pi^2}u^2.
\]
Our last task is to show that in the limit (\ref
{eqninclusionExclusion}) is equivalent to its $m=0$ term.
By iterations of the Hadamard--Fisher inequality, Lemma~\ref{lemFischer},
\[
\rho^{\U(n)}_{2k+m}(\theta_1,x_1,\ldots,\theta_k,x_k,y_1,\ldots
,y_m)\leq
\rho^{\U(n)}_{2k}(\theta_1,x_1,\ldots,\theta_k,x_k)
\prod_{i=1}^m\rho^{\U(n)}_1(y_i).
\]
The contribution of all terms with $m\geq1$ in (\ref
{eqninclusionExclusion}) is therefore bounded by
\[
\biggl(\int_{\theta_1+A_n}\dd x_1\cdots
\int_{\theta_k+A_n}\dd x_k\, \rho^{\U(n)}_{2k}(\theta_1,x_1,\ldots
,\theta
_k,x_k)\biggr)\sum_{m\geq1}\frac{1}{m!}\biggl(\int_{a_n}\rho^{\U
(n)}_1(y)\biggr)^m,
\]
where the integration domain $a_n=(\theta_1+A_n)\sqcup\cdots\sqcup
(\theta
_k+A_n)$ has size $O(n^{-4/3})$ and $\rho^{\U(n)}_1(y)=n$.
The first term of the above product converges, as previously proved (it
corresponds to $m=0$), so
the whole term goes to 0 as $n\to\infty$, concluding the proof.
\end{pf}
%
%
\begin{lem}[(Uniform boundness)]\label{lemuniformBounded}
There is a constant $c$ depending only on $A$ such that, for any $n\geq
1$ and
$(\theta_1,\ldots,\theta_k)\in\Omega^{\U(n)}$ [see (\ref{eqndefOmega})],
\[
\tilde\rho^{\U(n)}_k(\theta_1,\ldots,\theta_k)<c.
\]
\end{lem}
\begin{pf}
As previously mentioned, formula (\ref{eqninclusionExclusion})
is true whenever, for all distinct $i$ and $j$, $\theta_j\notin\theta
_i+A_n$, that is, $(\theta_1,\ldots,\theta_n)$ is in $\Omega^{\U(n)}$.
Using the Hadamard--Fisher inequality as in the proof of the previous
lemma, $\tilde\rho^{\U(n)}_k(\theta_1,\ldots,\theta_k)$ is therefore
bounded by
\[
\biggl(\int_{\theta_1+A_n}\dd x_1\cdots
\int_{\theta_k+A_n}\dd x_k\, \rho^{\U(n)}_{2k}(\theta_1,x_1,\ldots
,\theta
_k,x_k)\biggr)\sum_{m\geq0}\frac{1}{m!}\biggl(\int_{a_n}\rho^{\U
(n)}_1(y)\biggr)^m
\]
with $a_n=(\theta_1+A_n)\sqcup\cdots\sqcup(\theta_k+A_n)$.
Once again, the Hadamard--Fisher inequality gives
\[
\rho^{\U(n)}_{2k}(\theta_1,x_1,\ldots,\theta_k,x_k)\leq\prod
_{i=1}^k\rho
^{\U(n)}_2(\theta_i,x_i).
\]
This gives the upper bound, uniform in $(\theta_1,\ldots,\theta_n)\in
\Omega^{\U(n)}$,
\[
\biggl(\int_{A_n}\rho^{\U(n)}_2(0,x)\,\dd x\biggr)^k\sum_{m\geq
0}\frac
{1}{m!}\biggl(\int_{a_n}\rho^{\U(n)}_1(y)\biggr)^m
\]
converging to
\[
\biggl(\frac{1}{48\pi^2}\int_A u^2\,\dd u\biggr)^k
\]
as previously seen.
\end{pf}
\begin{rem*}
A better upper bound in the previous proof can be obtained as follows.
By a direct ensembles argument,
$\tilde\rho^{\U(n)}_k(\theta_1,\ldots,\theta_k)$ is bounded by
\[
\int_{\theta_1+A_n}\dd x_1\cdots
\int_{\theta_k+A_n}\dd x_k\, \rho^{\U(n)}_{2k}(\theta_1,x_1,\ldots
,\theta_k,x_k).
\]
This also comes from the fact that the inclusion--exclusion series (\ref
{eqninclusionExclusion}) is alternate.
We know by Fisher--Hadamard that this upper bound is lower than $(\int
_{A_n}\rho^{\U(n)}_2(0,x)\,\dd x)^k$, which is interpreted as follows:
$\tilde\rho^{\U(n)}_k$ converges to its limit from below, which is a
sign of repulsion before the asymptotic
independence.
\end{rem*}
%
%
\begin{lem}[(Negligible set)]\label{lemneglectSet}
Let $\overline{\Omega}{}^{\U(n)}$ be the complement of ${\Omega}^{\U(n)}$
in~$I^k$; see~(\ref{eqndefOmega}). Then
\[
\int_{\overline{\Omega}{}^{\U(n)}}\tilde\rho^{\U(n)}_k(\theta
_1,\ldots
,\theta_k)\,\dd\theta_1\cdots\dd\theta_k
\underset{n\to\infty}{\longrightarrow}0.
\]
\end{lem}
\begin{pf}
Let $(\theta_1,\ldots,\theta_k)\in\overline{\Omega
}{}^{\U
(n)}$, and note $\Theta=\{\theta_1,\ldots,\theta_k\}$ the set of
these points.
%
%
\begin{figure}[b]

\includegraphics{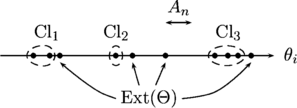}

\caption{Maximal clusters.}\label{fig3}
%
%
%
%
%
%
%
%
\end{figure}
For notational convenience, one can suppose $\theta_1<\cdots<\theta_k$.
A~set of points $\theta_{s}<\cdots<\theta_{t}$ is said to be a cluster
of $\Theta$ if,
for all $s\leq k\leq t$, $\theta_{k+1}\in\theta_k+A_n$.
Let $\Ext(\Theta)$ be the set of points which cannot be included in a
maximal cluster (see Figure~\ref{fig3})
\[
\Ext(\Theta)=\{\theta_i\dvtx1\leq i\leq k, (\theta_i+A_n)\cap\Theta
=\varnothing\}.
\]
This way we get a partition
\[
\Theta=\Ext(\Theta)\bigsqcup_{i=1}^\ell\Cl_i,
\]
where there are $\ell$ maximal clusters $\Cl_1,\ldots,\Cl_\ell$.
Suppose that $\Ext(\Theta)=\{\theta_{i_1}<\cdots<\theta_{i_p}\}$ where
$p=|{\Ext}(\Theta)|$.\vadjust{\goodbreak}
Then the following obvious bound holds:
\begin{eqnarray*}
\tilde\rho^{\U(n)}_k(\theta_1,\ldots,\theta_k)&\leq&\int_{\theta
_{i_1}+A_n}\dd x_1\cdots\int_{\theta_{i_p}+A_n}\dd x_p\,
\rho_{k+p}^{\U(n)}(\theta,\ldots,\theta_k,x_1,\ldots,x_p)\\
&\leq&
\prod_{j=1}^\ell\rho_{|{\Cl_i}|}^{\U(n)}(\Cl_i)\prod
_{j=1}^p\int_{\theta_{i_j}+A_n}
\rho^{\U(n)}_2(\theta_{i_j},x_j)\,\dd x_j,
\end{eqnarray*}
where we used the Hadamard--Fisher inequality. This last product of $p$
elements is
bounded, uniformly in $(\theta_1,\ldots,\theta_k)$, as it is equal to
the $p$th power of
\[
\int_{A_n}
\rho^{\U(n)}_2(0,x)\,\dd x,
\]
which converges.
Concerning the first product, we analyze each term in the following way.
For any $1\leq t\leq k-1$, the correlation $\rho^{\U(n)}_t$ is a
$t\times t$ determinant
of $K^{\U(n)}$.
Suppose that all arguments of $\rho^{\U(n)}_t$ are included in an
interval of size~$c_n$.
By subtracting the first row from the $t-1$ others, one obtains that
the first row is $O(n)$,
and the others are $O(n^2c_n)$ because $|K'^{\U(n)}|_\infty=O(n^2)$.
Hence $\rho^{\U(n)}_t=O(n^{2t-1}c_n^{t-1})$.
As a consequence, as all points of a cluster $\Cl_i$ are within
distance $O(n^{-4/3})$,
\[
\rho_{|{\Cl_i}|}^{\U(n)}(\Cl_i)=O\bigl(n^{
({2}/{3})|{\Cl
_i}|+{1/3}}\bigr).
\]
This leads to the upper bound
\[
\tilde\rho^{\U(n)}_k(\theta_1,\ldots,\theta_k)=
O\bigl(n^{({2}/{3})(k-p)+{\ell}/{3}}\bigr).
\]
But the size of the set of points $(\theta_1,\ldots,\theta_k)\in
\overline
{\Omega}{}^{\U(n)}$ with such a clusters configuration is
$O(n^{-({4}/{3})(k-p)})$, because there needs to be $k-p$ small gaps
[i.e., of size $O(n^{-4/3})$] between successive
$\theta_k$'s.
Hence the total contribution of such clusters to the integral over
$\overline{\Omega}{}^{\U(n)}$ is
\[
O\bigl(n^{{\ell}/{3}-({2}/{3})(k-p)}\bigr),
\]
which goes top 0 because $\ell\leq k-p$ (with equality if all clusters
have only one point) and
$k-p<2(k-p)$ [$k=p$ is not possible: this would mean that $(\theta
_1,\ldots,\theta_k)\notin\overline{\Omega}{}^{\U(n)}$].
\end{pf}

As previously explained, the three above lemmas complete the proof of
Theorem~\ref{thmunitary}.
\begin{pf*}{Proof of Corollary~\ref{corsmallestUnitary}} Note that
the events
$\{x_\ell<n^{4/3} \tau^{(n)}_\ell<y_\ell,1\leq\ell\leq k\}$ and
\begin{eqnarray*}
&&\bigl\{
\chi^{(n)}((72\pi)^{1/3}(x_k,y_k),(0,2\pi))\geq1,\\
&&\qquad\chi^{(n)}((72\pi)^{1/3}(x_\ell,y_\ell),(0,2\pi))=1,1\leq\ell
\leq
k-1,\\
&&\hspace*{12pt}\qquad\chi^{(n)}((72\pi)^{1/3}(y_{\ell-1},x_\ell),(0,2\pi))=0,1\leq\ell
\leq
k \bigr\}
\end{eqnarray*}
are almost surely the same ($y_0=0$).
The independence property of the limit Poisson point process $\chi$ in
disjoint subsets therefore yields
\begin{eqnarray*}
&&\Prob\bigl(x_\ell<n^{4/3} \tau^{(n)}_\ell<y_\ell,1\leq\ell
\leq
k\bigr)\\
&&\qquad\underset{n\to\infty}{\longrightarrow}
\bigl(1-e^{-(y_k^3-x_k^3)}\bigr)
\prod_{\ell=1}^{k-1}(y_\ell^3-x_\ell^3)e^{-(y_\ell^3-x_\ell^3)}
\prod_{\ell=1}^{k}e^{-(x_\ell^3-y_{\ell-1}^3)},
\end{eqnarray*}
where we noted that, for any interval $(a,b)$,
$\chi^{(n)}((72\pi)^{1/3}(a,b),(0,2\pi))$ is a Poisson random
variable with
parameter $(b^3-a^3)$. A straightforward simplification of the above
products gives the expected result.
Concerning the limiting density of $n^{4/3} \tau^{(n)}_k$, we proceed
in two steps. First, from formula (\ref{eqnjointLarge}),
the joint density of $n^{4/3}(\tau^{(n)}_1,\ldots,\tau^{(n)}_k)$ is, in
the limit and on the simplex $0<u_1<\cdots<u_k$, proportional to
\[
u_1^2u_2^2\cdots u_k^2e^{-u_k^3}.
\]
Consequently,
\begin{eqnarray*}
\Prob(\tau_k<x)&=&
c_k\int_0^x u_k^{2}e^{-u_k^3}\int_{0<u_1<\cdots
<u_{k-1}<u_k}u_1^2u_2^2\cdots u_{k-1}^2\,\dd u_1\cdots\dd u_{k-1}\\
&=&c_k\int_0^x u_k^{2+3(k-1)}e^{-u_k^3}\int_{0<v_1<\cdots
<v_{k-1}<1}v_1^2v_2^2\cdots v_{k-1}^2\,\dd v_1\cdots\dd v_{k-1},
\end{eqnarray*}
so the density of $\tau_k$ is proportional to $x^{3k-1}e^{-x^3}$.
\end{pf*}

\subsection{The Gaussian unitary ensemble}
The small gaps asymptotics for the $\Gue$ are obtained exactly in the
same way as for the unitary group. The only difference is that the
determinantal kernel
is not translation invariant anymore, leading to some complications.
More precisely, let $(h_n)$ be the Hermite polynomials, more precisely
the successive monic
orthogonal polynomials with respect to the Gaussian weight
$e^{-x^2/2}\,\dd x$. Following~\cite{AGZ}, where the following results on
the determinantal aspect of the $\Gue$ can be found,
we introduce the functions
\[
\psi_k(x)=\frac{e^{-x^2/4}}{\sqrt{\sqrt{2\pi}k!}} h_k(x).
\]
Then the set of points $\{\lambda_1,\ldots,\lambda_n\}$ with law
(\ref
{eqnGUE}) is a determinantal point process
with kernel (with respect to the Lebesgue measure on $\RR$) given by
%
%
\begin{eqnarray}\label{eqnkernelHermite}
K^{\Gue(n)}(x,y)&=&n\frac{\psi_n(x\sqrt{n})\psi_{n-1}(y\sqrt
{n})-\psi
_{n-1}(x\sqrt{n})\psi_n(y\sqrt{n})}{x-y}\nonumber\\
&=&n^{3/2}\biggl(\psi_{n-1}\bigl(y\sqrt{n}\bigr)\int_0^1\psi_n'\bigl(tx\sqrt
{n}+(1-t)y\sqrt{n}\bigr)\,\dd t\\
&&\hspace*{24pt}{}-\psi_{n}\bigl(y\sqrt{n}\bigr)\int
_0^1\psi
_{n-1}'\bigl(tx\sqrt{n}+(1-t)y\sqrt{n}\bigr)\,\dd t\biggr).
\nonumber
\end{eqnarray}
We will now discuss why Lemma~\ref{lemnoSucSmall} through Lemma \ref
{lemneglectSet} still hold in the case of the $\Gue$ ensemble,
restricted to $(-2+\epsilon_0,2-\epsilon_0)$ for a given $\epsilon_0>0$.
Note that the Hadamard--Fischer inequality, an important tool for the
proof of all the following lemmas, still holds for the Gaussian unitary
ensemble because the set of its eigenvalues is a determinantal process.
\begin{itemize}
\item We now note
\[
\xi^{(n)}=\sum_{i=1}^n\delta_{\lambda_i}\one_{|\lambda
_i|<2-\epsilon_0}
\]
and its thinning $\tilde\xi^{(n)}$
obtained from $\xi^{(n)}$ by only keeping the eigenvalues $\lambda_k$
for which $\xi^{(n)}(\lambda_k+n^{-4/3}A=1)$. The absence of successive
small gaps, Lemma~\ref{lemnoSucSmall}, only requires in its proof that
$|K^{\Gue(n)}|_\infty=\OO(n)$, $|\partial_x K^{\Gue(n)}|_\infty
=\OO
(n^2)$. This is proved in Lemma~\ref{leminfiniteNorms}.
\item The analog of the simple convergence of the correlation function
associated to the thinned point process, Lemma \ref
{lemsimpleConvergence}, is now, for any distinct points $\lambda
_1,\ldots,\lambda_k$
in $(-2+\epsilon_0,2-\epsilon_0)$,
\[
\tilde\rho^{\Gue(n)}_k(\lambda_1,\ldots,\lambda_k)\underset{n\to
\infty
}{\longrightarrow}\biggl(\frac{1}{48\pi^2}\int_A u^2\,\dd u
\biggr)^k\prod
_{i=1}^k(4-\lambda_i^2)^2.
\]
The proof, in the same way as the unitary case, requires uniform bounds
on the partial derivatives of the kernel, given in Lemma \ref
{leminfiniteNorms}, and the asymptotics
\[
\int_{\lambda_i+A_n}\rho^{\Gue(n)}_2(\lambda_i,x)\,\dd x\underset
{n\to
\infty}{\longrightarrow}\biggl(\frac{1}{48\pi^2}\int_Au^2\,\dd
u
\biggr)(4-\lambda_i^2)^2,
\]
which is a direct consequence of Lemma~\ref{lemTaylorExpansion}.
\item The uniform boundness result, Lemma~\ref{lemuniformBounded},
only requires the Hadamard--Fischer inequality, Lemma \ref
{lemFischer}, which applies to any determinantal point process.
\item Finally, we give the analog of Lemma~\ref{lemneglectSet} in the
following way. Let
\[
\Omega^{\Gue(n)}=\{(\lambda_1,\ldots,\lambda_k)\in I^k\dvtx\lambda
_i\notin
\lambda_j+A_n,1\leq i,j\leq k\},
\]
and $\overline{\Omega}{}^{\Gue(n)}$ be its complement in $I^k$, where $I$
is any Borel set in $(-2+\epsilon_0$, $2-\epsilon_0)$.
\[
\int_{\overline{\Omega}{}^{\Gue(n)}}\tilde\rho^{\Gue(n)}_k(\lambda
_1,\ldots
,\lambda_k)\,\dd\lambda_1\cdots\dd\lambda_k
\underset{n\to\infty}{\longrightarrow}0.
\]
The only estimate necessary for this result is that the first-order
partial derivative of the kernel is uniformly $\OO(n)$, which is one of
the estimates of the following Lemma
\ref{leminfiniteNorms}.
\end{itemize}
%
%
\begin{lem}\label{leminfiniteNorms}
Let $\epsilon_0>0$. Uniformly for
$x,y\in(-2+\epsilon_0,2-\epsilon_0)$,\break
$K^{\Gue(n)}(x$, $y)$ is $\OO(n)$, the first-order partial derivatives of
$K^{\Gue(n)}$ are $\OO(n^2)$ and
the second-order ones are $\OO(n^3)$.

Moreover, under the additional condition $|x-y|>\delta>0$,
$K^{\Gue(n)}(x,y)$ is uniformly bounded, by a constant depending only
on $\epsilon_0$ and $\delta$.
\end{lem}
\begin{pf}
From the Plancherel--Rotach asymptotics for the Hermite polynomials
(Theorem 8.22.9 in~\cite{Szego}),
for any nonnegative integer $k$, $\psi_{n-k}(\sqrt{n}x)$ is $\OO
(1/\sqrt{n})$, uniformly in $x\in(-2+\epsilon_0,2-\epsilon_0)$.
Consequently, if $|x-y|>\delta$, from the first line of (\ref
{eqnkernelHermite}), $K^{\Gue(n)}$ is uniformly $\OO(1)$.

To prove the first assertion of the lemma, we use the stability
property of the functions $(\psi_k)$ by derivation (see~\cite{AGZ},
Lemma 3.2.7),
%
%
\begin{equation}\label{eqnstability}
\psi_n'(x)=-\frac{x}{2}\psi_n(x)+\sqrt{n}\psi_{n-1}(x).
\end{equation}
Injecting this expression of $\psi_n'$ in the last two lines of (\ref
{eqnkernelHermite}) and using $\psi_n=\OO(1/\sqrt{n})$
yields $K^{\Gue(n)}(x,y) =\OO(n)$ uniformly for $x,y\in(-2+\epsilon
_0,2-\epsilon_0)$.
Iterating this procedure, any partial derivative of $K^{\Gue(n)}$ of
order $k$ are $n^{k+2}$ times a linear combination of factors of type
$\psi_{n-i}\int_0^1\psi_{n-j}$,
which shows that any partial derivative of order $k$ is $\OO(n^{k+1})$.
\end{pf}

The Plancherel--Rotach asymptotics for Hermite polynomials also yield a
precise evaluation of the correlation functions evaluated at close
points, in particular the
following result.
%
%
\begin{lem}\label{lemTaylorExpansion}
Let $\epsilon_0>0$, $c>0$. Then as $n\to\infty$, uniformly in $x\in
(-2+\epsilon_0$, $2-\epsilon_0)$, $|u|<c n^{-4/3}$,
\[
\rho^{\Gue(n)}_2(x,x+u)=\frac{1}{48\pi^2} n^4 (4-x^2)^2u^2+\OO
(n).
\]
\end{lem}
\begin{pf}
The intuition for this result is as follows. From the well-known
convergence to the sine kernel,
\[
\frac{1}{n \rho_{\mathrm{sc}}(x)} K^{\Gue(n)}\biggl(x,x+\frac{v}{n\rho
_{\mathrm{sc}}(x)}\biggr) \underset{n\to\infty}{\longrightarrow} \frac
{\sin{(\pi
v)}}{\pi v}.
\]
If this convergence is sufficiently uniform, one expects for small $u$
a Taylor expansion
\[
K^{\Gue(n)}(x,x+u)\approx\frac{\sin{(\pi n\rho_{\mathrm{sc}}(x)
u)}}{\pi u}\approx n\rho_{\mathrm{sc}}(x)-\frac{\pi^2}{6}\rho_{\mathrm{sc}}(x)^3n^3 u^2.
\]
Injecting this expansion in
\[
\rho^{\Gue(n)}_2(x,x+u)=K^{\Gue(n)}(x,x)K^{\Gue
(n)}
(x+u,x+u)-K^{\Gue(n)}(x,x+u)^2
\]
gives the expected result. To justify the accuracy of the sine kernel
approximation for $u$ close to 0,
we rely on Corollary 1 in~\cite{DelYao} (more general
asymptotics in the full complex plane for, in particular, Hermite
polynomials were
obtained in~\cite{DKMVZ}):
using the Plancherel--Rotach asymptotics of Hermite polynomials the
authors show that\footnote{Note that the normalization in~\cite{DelYao}
is different, as their semicircle law is supported on
$(-\sqrt{n},\sqrt{n})$.}
\[
K^{\Gue(n)}(x,y)=\frac{1}{2\pi\sqrt{\sin\omega\sin
\theta
}}\frac{\sin(n(a_\theta-a_\omega))}{\sin(\omega
-\theta
)}+\OO(1)
\]
uniformly for $x,y$ in $(-2+\epsilon_0, 2-\epsilon_0)$
where $x=2\cos\omega$, $y=2\cos\theta$, $0\leq\theta$, $\omega\leq
\pi$,
$a_\omega=\sin(2\omega)-2\omega$, $a_\theta=\sin(2\theta
)-2\theta$. As
$a_\theta-a_\omega=2(\omega-\theta)\sin^2\theta+\OO((\omega
-\theta)^2)$ and
$\omega-\theta=(y-x)/\sqrt{4-x^2}+\OO((x-y)^2)$, a simple expansion yields
%
%
\begin{eqnarray}\label{eqnGUEexpansion}
&&\frac{1}{\rho_{\mathrm{sc}}(x)}K^{\Gue(n)}(x,y)\nonumber\\
&&\qquad=\frac{1}{2\pi^2\rho_{\mathrm{sc}}(x)
\sqrt{\rho_{\mathrm{sc}}(x)\rho_{\mathrm{sc}}(y)}}\nonumber\\[-8pt]\\[-8pt]
&&\qquad\quad{}\times\frac
{\sin
(2n(\omega-\theta)\sin^2\theta+\OO(n(\omega-\theta
)^2))}{\sin
(({x-y})/{\sqrt{4-y^2}}+\OO((x-y)^2))}+\OO
(1)\nonumber\\
&&\qquad=\frac{\sin(n\pi(x-y)\rho_{\mathrm{sc}}(x)+\OO(n(x-y)^2))}{\pi
(x-y)\rho_{\mathrm{sc}}(x)+\OO((x-y)^2)}+\OO(1).
\nonumber
\end{eqnarray}
For $x-y=\OO(n^{-4/3})$, this implies
\[
\frac{1}{\rho_{\mathrm{sc}}(x)}K^{\Gue(n)}(x,y)=n-\frac
{1}{6}n^3\pi
^2(x-y)^2\rho_{\mathrm{sc}}(x)^3+\OO(1)
\]
and the expected result for $\rho^{\Gue(n)}_2(x,x+u)$.
\end{pf}

The above discussion completes the proof of Theorem~\ref{thmG}. Its
Corollary~\ref{corsmallestGUE}
follows exactly in the same way as Corollary~\ref{corsmallestUnitary}
from Theorem~\ref{thmunitary}.

Corollary~\ref{corsmallestlocation} is a straightforward consequence
of the scission of the limiting L\'evy measure:
if $\tilde\chi=\{(a_j,b_j)\}$ is a Poisson point process with measure
$\E\tilde\chi(A\times B)=\mu_1(A)\mu_2(B)$, conditionally to the
$b_k$'s the $a_k$'s are distributed independently of
each other, independently of the $b_k$'s, as a Poisson point process
with intensity proportional to $\mu_1$. In particular, in our
situation, the abscissa associated to the minimal ordinate
is distributed with density proportional to
\[
(4-x^2)^2.
\]

\section{Large gaps}\label{sec3}

\subsection{The unitary group: Asymptotics of Toeplitz
determinants}\label{subsecRH}
To evaluate the extreme gaps, we first investigate the queuing
distribution of one given spacing. A large part of the literature
concerns the probability of having no
eigenvalues in a given interval (e.g.,~\cite{Dyson1,Widom}). A rigorous
derivation of the queuing distribution
for the large gaps was only given recently, thanks to the steepest
descent method for Riemann--Hilbert problems
\cite{Deift,Kra} or by operator theory tools for Toeplitz determinants
\cite{Ehr}.

The link with the nearest neighbor distribution is given by the
following lemma, explained
in~\cite{Mehta}, Appendix A.8.
%
%
\begin{lem}\label{lemMeh}
Let $u\sim\mu_{\U(n)}$, with eigenangles $0\leq\theta_1<\cdots
<\theta
_n<2\pi$, and $\xi^{(n)}=\sum_{i=1}^n\delta_{\theta_i}$. Then
\[
n\Prob(\theta_2-\theta_1>x)=-\frac{\dd}{\dd x}\Prob
\bigl(\xi
^{(n)}(0,x)=0\bigr).
\]
\end{lem}

Moreover, the probability of having no eigenvalues in an arc of size
$2\alpha$ is equal to the Toeplitz determinant (this is a consequence
of Heine's formula)
\[
D_n(\alpha)=\underset{1\leq j,k\leq n}{\det}\biggl(\frac{1}{2\pi
}\int
_\alpha^{2\pi-\alpha}e^{\ii(j-k)\theta}\,\dd\theta\biggr).
\]

All the asymptotics we need in the following are direct consequences of
the precise analysis of $D_n(\alpha)$ given by Deift et al.~\cite{Deift}
and Krasovsky~\cite{Kra}.
More precisely they prove that for some sufficiently large $s_0$ and
any $\epsilon>0$, uniformly in
$s_0/n<\alpha<\pi-\epsilon$,
%
%
\begin{eqnarray}
\label{eqnest1}
\log D_n(\alpha)&=&n^2\log\cos\frac{\alpha}{2}-\frac{1}{4}\log
\biggl(n\sin\frac{\alpha}{2}\biggr)+c_0+\OO\biggl(\frac{1}{n\sin(\alpha
/2)}\biggr),\\
\label{eqnest2}\qquad
\frac{\dd}{\dd\alpha}\log D_n(\alpha)&=&-\frac{n^2}{2}\tan\frac
{\alpha
}{2}-\frac{1}{8}\cot\frac{\alpha}{2}+\OO\biggl(\frac{1}{n\sin
^2(\alpha
/2)}\biggr)
\end{eqnarray}
for an explicit constant $c_0$, which remained conjectural for 20 years.
From this, we can give upper bounds on
$n\Prob(\theta_2-\theta_1>u)$, the expectation of the number
of gaps greater than $u$.
%
%
\begin{lem}\label{lemboundExpectation}
Given any $a>0$ and $\epsilon>0$, uniformly in $u\in(a\sqrt{\log
n}/n$, $2\pi-\epsilon)$, as $n\to\infty$
\[
n\Prob(\theta_2-\theta_1>u)\leq u n^{2+\o(1)}
e^{-n^2
{u^2}/{32}}.
\]
\end{lem}
\begin{pf}
From Lemma~\ref{lemMeh},
\[
n\Prob(\theta_2-\theta_1>u)=-\frac{1}{2}\frac{\dd
}{\dd\alpha
}D_n(\alpha)=-\frac{1}{2}\biggl(\frac{\dd}{\dd\alpha}\log
D_n(\alpha
)\biggr) D_n(\alpha)
\]
evaluated for $u=2\alpha$.
The first term is evaluated thanks to (\ref{eqnest2}): it is $\OO(
u
n^{2+\o(1)})$ uniformly in $(a\sqrt{\log n}/n,\pi-\epsilon)$.
Concerning the second term, (\ref{eqnest1}) and the inequality $\log
\cos x\leq-x^2/2$ on $[0,\pi/2)$ imply that it is\break $e^{-n^2\alpha
^2/8}n^{\o(1)}$,
completing the proof.
\end{pf}

Note that substituting $u=2\alpha=\frac{\sqrt{\lambda\log n}}{n}$ in
the asymptotics (\ref{eqnest1}) and (\ref{eqnest2}) yields
%
%
\begin{equation}\label{eqntransition}
n\Prob(\theta_2-\theta_1>u)=n^{1-{\lambda
}/{32}+\o(1)},
\end{equation}
where the $\o(1)$ terms are uniform in $\lambda\in[a,b]$, for any given
positive $a$ and~$b$.
Hence, if $\lambda>32$, the expected number of gaps greater than
$\frac
{\sqrt{\lambda\log n}}{n}$ goes to $0$. If
$\lambda<32$, this goes to $\infty$.
The transition in the number of large gaps at $\lambda=32$ is a strong
clue for Theorem~\ref{thmlargestGaps},
but it is not sufficient: concluding the proof requires additional
knowledge about the correlations between gaps, Lemma~\ref{lemNegCor},
assumed for the moment.
Take $p>0$. We denote
\[
X_n=\frac{n}{\sqrt{32\log n}} \mathcal{T}^{(n)}_{\ell_n},
\]
where $\ell_n=n^{\o(1)}$. To prove that $X_n-1$ converges to 0 in $L^p$,
we pick some arbitrarily small $\epsilon>0$, bound
$\E(|X_n-1|^p\one_{1-\epsilon<X_n<1+\epsilon})$
by $\epsilon^p$, and we need to prove that both of the following terms
converge to 0:
%
%
\begin{eqnarray}\label{eqn3terms}
\E(|X_n-1|^p\one_{X_n<1-\epsilon})&\leq&\Prob
(X_n<1-\epsilon
),\nonumber\\
\E(|X_n-1|^p\one_{X_n>1+\epsilon}) &\leq& \E
\bigl(|X_n-1|^p\one_{X_n>1+\epsilon,\mathcal{T}^{(n)}_{1}<\pi/2}\bigr)\\
&&{}+\biggl(\frac{2\pi n}{\sqrt{32\log n}}\biggr)^p\Prob
\bigl(\mathcal
{T}^{(n)}_{1}>\pi/2\bigr)\nonumber
\end{eqnarray}
(here the value $\pi/2$ is arbitrary, any angle strictly smaller than
$\pi$ would be appropriate for this proof).
First, decompose the unit circle into 8 fixed angles of size $\pi/4$.
If $\mathcal{T}^{(n)}_{1}>\pi/2$, one of these arcs is free of eigenvalues.
From (\ref{eqnest1}), the probability of this event decreases
exponentially, so (\ref{eqn3terms}) converges to 0.

Moreover, integrating by parts,
\begin{eqnarray*}
&&\E\bigl(|X_n-1|^p\one_{X_n>1+\epsilon,\mathcal{T}^{(n)}_{1}<\pi
/2}
\bigr)\\
&&\qquad=
\int_{1+\epsilon}^{\infty} p(u-1)^{p-1}\Prob\bigl(X_n>u,\mathcal
{T}^{(n)}_{1}<\pi/2\bigr)\,\dd u\\
&&\qquad\quad{}
+p \epsilon^{p-1}\Prob\bigl(X_n>1+\epsilon,\mathcal{T}^{(n)}_{1}<\pi
/2\bigr)\\
&&\qquad\leq\int_{1+\epsilon}^{({n}/{\sqrt{32\log n}})({\pi}/{2})}
p(u-1)^{p-1}\Prob(X_n>u)\,\dd u
+p \epsilon^{p-1}\Prob(X_n>1+\epsilon),
\end{eqnarray*}
because $X_n$ needs to be shorter than $\frac{n}{\sqrt{32\log
n}}\mathcal{T}^{(n)}_{1}$.
The probability that $X_n$ is greater than $u$ is obviously shorter
than $n\Prob(\theta_2-\theta_1>u\sqrt{32 \log n}/n)$,
the expectation of the number of gaps greater than $u\sqrt{32 \log
n}/n$. Hence
the above quantity goes to $0$ thanks to the uniform estimate of Lemma
\ref{lemboundExpectation}.

Finally, showing that $\Prob(X_n<1-\epsilon)\to0$ requires an
additional argument, the negative correlation property for empty sets
events, Lemma~\ref{lemNegCor}.
In particular, this negative correlation implies the following result.
%
%
\begin{lem}\label{lemVarianceExpectation}
Consider a set of disjoint arcs on the unit circle. Let $M_n$ be the
number of such intervals free of eigenangles, that is, those $I_k$'s
such that
\mbox{$
\xi^{(n)}(I_k)=0
$}.
Then $\Var(M_n)\leq\E(M_n)$.
\end{lem}
\begin{pf}
This is a straightforward consequence of Lemma~\ref{lemNegCor}: as
$\xi
^{(n)}$ is a determinantal point process,
for disjoints $I$ and $J$, $\Prob(\xi^{(n)}(I\cup J)=0)\leq\Prob
(\xi
^{(n)}(I_j)=0)\Prob(\xi^{(n)}(I_k)=0)$. Hence, noting $m_n$ the number
of initial arcs,
\begin{eqnarray*}
\E(M_n^2)&=&\sum_{1\leq j,k\leq m_n}\Prob\bigl(\xi^{(n)}(I_j\cup
I_k)=0\bigr)\\
&\leq&2 \sum_{1\leq j<k\leq m_n}\Prob\bigl(\xi^{(n)}(I_j)=0\bigr)\Prob\bigl(\xi
^{(n)}(I_k)=0\bigr)+\sum_{1\leq j\leq m_n}\Prob\bigl(\xi^{(n)}(I_j)=0\bigr)\\
&\leq& \E(M_n)^2+\E(M_n)
\end{eqnarray*}
as expected.
\end{pf}

Consider now a number $m_n$ of disjoint intervals $I_1,\ldots,I_{m_n}$
of length $(1-\epsilon)\frac{\sqrt{32\log n}}{n}$ in $(0,2\pi)$.
We can find $m_n=\lfloor2\pi n/\sqrt{32\log n}\rfloor$ of them. Let
$M_n$ be the number of such intervals free of eigenangles,
that is, those $I_k$'s such that
$
\xi^{(n)}(I_k)=0
$.
If there are less than $\ell_n$ gaps larger than $(1-\epsilon)\frac
{\sqrt{32\log n}}{n}$, either there are less than
$\ell_n$ intervals $I_k$'s free of eigenangles, or there is a gap
between successive eigenangles containing two intervals,
\[
\Prob(X_n<1-\epsilon)\leq\Prob(M_n<\ell_n)+\Prob\biggl(\mathcal
{T}^{(n)}_1\geq2(1-\epsilon)\frac{\sqrt{32\log n}}{n}\biggr).
\]
This last term is bounded by
the expectation of the number of gaps greater than $2(1-\epsilon)\frac
{\sqrt{32\log n}}{n}$; from (\ref{eqntransition})
this goes to 0 if $2(1-\epsilon)>1$ (true for $\epsilon$ sufficiently small).

Concerning $\Prob(M_n<\ell_n)$, first note that by the estimate (\ref
{eqntransition}), $\E(M_n)=n^{1-(1-\epsilon)^2+\o(1)}$, so $\ell
_n=\o
(\E(M_n))$. This\vadjust{\goodbreak}
allows to use Chebyshev's inequality, for sufficiently large $n$.
\begin{eqnarray*}
\Prob(M_n<\ell_n)&\leq&\Prob\bigl(|M_n-\E(M_n)|>\E(M_n)-\ell_n\bigr)\\
&\leq&
\frac{\Var
(M_n)}{(\E(M_n)-\ell_n)^2}\leq\frac{\E(M_n)}{(\E(M_n)-\ell_n)^2},
\end{eqnarray*}
where we used Lemma~\ref{lemVarianceExpectation} in the last inequality.
This last term is equivalent to $1/\E(M_n)$, thus going to $0$, which
completes the proof.

\subsection{The GUE: Comparison of Fredholm determinants}
For the proof of Theorem~\ref{thmlargestGapsGUE}, precise asymptotics
like (\ref{eqnest1}) and (\ref{eqnest2})
related to unitary groups are not available in the $\Gue$ context. This
difficulty can be overcome, our main observation being Lemma
\ref{lemcompFredholm}: the probability that an interval is free of
eigenvalues is equivalent in the $\Gue(n)$ and $\U(n)$ cases, up to a
normalization, if the interval size is shorter than the expected
extreme gap size.
The proof relies on a comparison of the Fredholm determinants
associated to $K^{\Gue(n)}$ and $K^{\U(n)}$.

The rest of the proof is similar to the one concerning $\U(n)$. Indeed,
consider an $\epsilon>0$, $p>0$, $\mathcal{\tilde T}^{(n)}_{\ell_n}$
the $\ell_n$th largest gap in $I$, and
\[
\tilde X_n=\frac{\mathcal{\tilde T}^{(n)}_{\ell_n}}{t_n},\qquad t_n=\frac
{\sqrt{32\log n}}{n \inf_I\sqrt{4-x^2}}.
\]
We then decompose
\[
\E(|\tilde X_n-1|^p)\leq\epsilon^p+\Prob(\tilde
X_n<1-\epsilon)+\E(|\tilde X_n-1|^p\one_{\tilde X_n>1+\epsilon
}).
\]
This last expectation is, by integration by parts,
\[
\int_{\epsilon}^\infty p v^{p-1}\Prob(\tilde X_n>1+v)\,\dd v+p\epsilon
^{p-1}\Prob(\tilde X_n>1+\epsilon).
\]
The probability $\Prob(\tilde X_n>1+v)$ is lower than the expectation
of the number of gaps greater than $(1+v)t_n$. Lemma~\ref{leminTheTailGUE}
therefore yields $\E(|\tilde X_n-1|^p\one_{\tilde
X_n>1+\epsilon
})\to0$ as $n\to\infty$.

Concerning $\Prob(\tilde X_n<1-\epsilon)$, we proceed as for the
unitary group: for $I=(a,b)$ with $a<b$, consider
$\tilde m_n=\lfloor(b-a)/((1-\epsilon) t_n)\rfloor$ disjoint intervals
of length $(1-\epsilon)t_n$ included in $I$.
Let $\tilde M_n$ be the number of these intervals free of eigenvalues.
Then
\[
\Prob(\tilde X_n<1-\epsilon)\leq\Prob(\tilde M_n<\ell_n)+\Prob
\bigl(\mathcal
{\tilde T}^{(n)}_{1}>2(1-\epsilon)t_n\bigr).
\]
From Lemma~\ref{leminTheTailGUE} this last probability goes to $0$ as
$n\to\infty$ if $2(1-\epsilon)>1$, $\epsilon<1/2$.
Moreover, from Lemma~\ref{lemexpectGUE},
%
%
\begin{equation}\label{eqnexpboundGUE}
\E(\tilde M_n)\gg n^\delta
\end{equation}
for some $\delta>0$ depending only on $\epsilon$ and $I$. Moreover,
from the negative correlation property Lemma~\ref{lemNegCor} and
the same reasoning as Lemma~\ref{lemVarianceExpectation},
%
%
\begin{equation}\label{eqnVarianceExpectationGUE}
\Var(\tilde M_n)\leq\E(\tilde M_n).\vadjust{\goodbreak}
\end{equation}
As $\ell_n=n^{\o(1)}$, from (\ref{eqnexpboundGUE}) $\E(\tilde
M_n)-\ell
_n>0$ for sufficiently large $n$, which allows us to use Chebyshev's inequality,
\begin{eqnarray*}
\Prob(\tilde M_n<\ell_n)&\leq&\Prob\bigl(|\tilde M_n-\E(\tilde M_n)|>\E
(\tilde M_n)-\ell_n\bigr)\\
&\leq&\frac{\Var(\tilde M_n)}{(\E(\tilde
M_n)-\ell
_n)^2}\leq\frac{\E(\tilde M_n)}{(\E(\tilde M_n)-\ell_n)^2},
\end{eqnarray*}
the last inequality being (\ref{eqnVarianceExpectationGUE}).
From (\ref{eqnexpboundGUE}) this last term is equivalent to $1/\E
(\tilde M_n)$ and going to 0, as $n\to\infty$, completing the proof.
%
%
\begin{lem}\label{lemasymptSine}
Let $\delta_n=\o(1)$. The following asymptotics hold for the unitary
group and $\Gue$ kernels:
\begin{longlist}[(2)]
\item[(1)] Uniformly for $x,y$ in $(0,2\pi)$ and $|x-y|=\OO(\delta_n)$,
\[
\frac{2\pi}{n}K^{\U(n)}(x,y)-\frac{\sin(n
({x-y})/{2}
)}{n({x-y})/{2}}=\OO\biggl(\frac{\delta_n}{n}\biggr).
\]
\item[(2)] Let $\epsilon_0>0$. Uniformly for $x,y$ in $(-2+\epsilon
_0,2-\epsilon_0)$ and $|x-y|=\OO(\delta_n)$,
\[
\frac{1}{n\rho_{\mathrm{sc}}(x)}K^{\Gue(n)}(x,y)-\frac{\sin(n\pi\rho
_{\mathrm{sc}}(x)(x-y))}{n\pi\rho_{\mathrm{sc}}(x)(x-y)}=\OO\biggl(\frac
{1}{n}
\biggr)+\OO(\delta_n)+\OO(n\delta_n^2).
\]
\end{longlist}
\end{lem}
\begin{pf}
For the unitary group, the kernel is explicit so
\begin{eqnarray*}
&&\frac{2\pi}{n}K^{\U(n)}(x,y)-\frac{\sin(n
({x-y})/{2}
)}{n({x-y})/{2}}\\
&&\qquad=
\sin\biggl(n\frac{x-y}{2}\biggr)\frac{1}{n({x-y})/{2}}
\biggl(\frac{
(({x-y})/{2})}{\sin(({x-y})/{2})}-1\biggr).
\end{eqnarray*}
As $|x-y|=\OO(\delta_n)\to0$, by expansion of $\sin$ at third order
$(\frac{x-y}{2})/\sin(\frac{x-y}{2})-1=\OO
((x-y)^2)$,
which completes the proof.

Concerning the Gaussian unitary ensemble, the same type of asymptotics
hold, being a direct consequence of formula (\ref{eqnGUEexpansion}).
Note that when taking $\delta_n=\OO(1/n)$,
the speed of convergence to the sine kernel is $1/n^2$ for the unitary
group, much better that $1/n$ for the $\Gue$, whose correlation kernel
is not translation invariant.
 \end{pf}
%
%
\begin{lem}\label{lemcompFredholm}
Let $\delta_n=\OO(\sqrt{\log n}/n)$, $\epsilon_0>0$. Then uniformly for
$x$ in $(-2+\epsilon_0$, $2-\epsilon_0)$,
\begin{eqnarray*}
&&\biggl|\Prob^{\Gue(n)}\biggl(\lambda_i\notin\biggl[x,x+\frac{\delta
_n}{\rho
_{\mathrm{sc}}(x)}\biggr],1\leq i\leq n\biggr)\\
&&\hspace*{38pt}{}-\Prob^{\U(n)}(\theta_i\notin[0,2\pi\delta_n
],1\leq
i\leq n)\biggr|\leq n^{\o(1)-1}.
\end{eqnarray*}
\end{lem}
\begin{pf}
By inclusion--exclusion, the probability that a determinantal point
process with kernel $K$ has no points in a measurable subset $A$ is the
Fredholm determinant (see, e.g., Lemma 3.2.4 in~\cite{AGZ})
%
%
\begin{equation}\label{eqnfredholm}
\det(\Id-K_A)=1+\sum_{k=1}^\infty\frac{(-1)^k}{k!}\int
_{A^k}\underset
{k\times k}{\det}(K(x_i,x_j))\,\dd x_1\cdots\dd x_k.
\end{equation}
%
To compare empty sets probabilities, we therefore need to compare
Fredholm determinants. A classical inequality (see, e.g.,~\cite{GGK},
Chapter IV, (5.14)) is
%
%
\begin{equation}\label{eqnineqTrace}
|{\det}(\Id+A)-\det(\Id+B)|\leq|A-B|_1e^{1+|A|_1+|B|_1},
\end{equation}
where $|T|_1$ is the trace norm of a nuclear operator $T$. However, for
positive (like $A$ or $B$)
operators $T(f)(x)=\int K(x,y)f(y)\,\dd y$ the trace norm is $\int
|K(x,x)|\,\dd x$,
but for nonpositive operators, like $A-B$, the trace norm is difficult
to evaluate, even for a compactly supported continuous kernel $K$.
However, in such a case, the Hilbert--Schmidt norm is computable.
%
%
\begin{equation}\label{eqnHSnorm}
|T|_2^2=\iint|K(x,y)|^2\,\dd x\,\dd y.
\end{equation}
This is the reason why, instead of the Fredhom determinant and the
inequality (\ref{eqnineqTrace}), we will use the modified
Carleman--Fredholm determinant
\[
\mathrm{det}_2(\Id+T)=\det(\Id+T)e^{-\Tr T}
\]
and the inequality
%
%
\begin{equation}\label{eqnineqHS}
|\mathrm{det}_2(\Id+A)-\mathrm{det}_2(\Id+B)|\leq|A-B|_2e^{
(|A|_2+|B|_2+1)^2/{2}},
\end{equation}
which can be found in~\cite{GGK}, Chapter IV, (7.11). For our purpose,
note that
from (\ref{eqnfredholm}), after a simple change of variables, the
probability that there are no eigenvalues in $[x,x+\delta_n/(\rho
_{\mathrm{sc}}(x))]$ is equal to
\[
1+\sum_{k=1}^\infty\frac{(-1)^k}{k!}\int_{(n\delta_n)^k}\underset
{k\times k}{\det}\biggl(\frac{1}{n\rho_{\mathrm{sc}}}K^{\Gue(n)}
\biggl(x+\frac
{y_i}{n\rho_{\mathrm{sc}}(x)},x+\frac{y_j}{n\rho_{\mathrm{sc}}(x)}\biggr)\biggr)\,\dd
y_1\cdots\dd y_k.
\]
In the same way, the probability that there are no eigenvalues in
$[0,2\pi\delta_n]$ is
\[
1+\sum_{k=1}^\infty\frac{(-1)^k}{k!}\int_{(n\delta_n)^k}\underset
{k\times k}{\det}\biggl(\frac{2\pi}{n}K^{\U(n)}\biggl(\frac{2\pi
}{n}y_i,\frac{2\pi}{n}y_j\biggr)\biggr)\,\dd y_1\cdots\dd y_k.
\]
Hence inequality (\ref{eqnineqHS}) will be applied with $A$ and $B$
integral operators with respective kernel
\[
A(u,v)=-\frac{1}{n\rho_{\mathrm{sc}}}K^{\Gue(n)}_{(0,n\delta_n)}
\biggl(x+\frac
{u}{n\rho_{\mathrm{sc}}(x)},x+\frac{v}{n\rho_{\mathrm{sc}}(x)}\biggr)
\]
and
\[
B(u,v)=-\frac{2\pi}{n}K^{\U(n)}_{(0,n\delta_n)}\biggl(\frac{2\pi
}{n}u,\frac{2\pi}{n}v\biggr).
\]
From Lemma~\ref{lemasymptSine}, the infinite norm between the two
kernels above is $\OO(n\delta_n^2)=\OO(\log n/n)$, so by (\ref
{eqnHSnorm}), integrating on a domain of area $(n\delta_n)^2=\OO
(\log n)$,
\[
|A-B|_2=\bigl(\OO\bigl((\log n)^2/n^2\bigr)\OO(\log n)\bigr)^{1/2}=\OO\bigl((\log n)^{3/2}/n\bigr).
\]
Moreover, consider a parameter $\alpha_n>0$, and decompose
\[
|A|_2^2=\iint|A(x,y)|^2\one_{|x-y|<\alpha_n}\,\dd x\,\dd y+\iint
|A(x,y)|^2\one_{|x-y|>\alpha_n}\,\dd x\,\dd y.
\]
From Lemma~\ref{lemasymptSine}, if $|x-y|>\alpha_n$, then $|A(x,y)|$
is smaller than $1/(\pi\alpha_n)+\OO(\log n/n)$, and when
$|x-y|<\alpha_n$,
it is bounded by $1+\OO(\log n/n)$. Consequently,
\[
|A|_2^2=\OO\bigl(\alpha_n\sqrt{\log n}\bigr)+\OO(\log n/\alpha_n^2)=\OO
((\log n)^{2/3})
\]
by choosing $\alpha_n=(\log n)^{1/6}$. In the same way, $|B|_2^2=\OO
((\log n)^{2/3})$. Hence
%
%
\begin{equation} \label{eqnineqDet2}
|\mathrm{det}_2(\Id+A)-\mathrm{det}_2(\Id+B)|= \OO\biggl(\frac{(\log
n)^{3/2}}{n}\biggr)e^{\OO((\log n)^{2/3})}\leq n^{\o(1)-1}.\hspace*{-30pt}
\end{equation}
Finally, using once again Lemma~\ref{lemasymptSine}, here on the diagonal,
%
%
\begin{equation}\label{eqnestimateTrace}
\Tr A=-\int A(x,x)\,\dd x=-n\delta_n+\OO\bigl((\log n)^{3/2}/n\bigr),
\end{equation}
and the same for $\Tr B$. Finally, write
\begin{eqnarray*}
&&|{\det}(\Id+A)-\det(\Id+B)|\\
&&\qquad\leq e^{\Tr A}|\mathrm{det}_2(\Id
+A)-\mathrm{det}_2(\Id+B)|+|e^{\Tr A-\Tr B}-1| |{\det}(\Id+B)|.
\end{eqnarray*}
Formulas (\ref{eqnineqDet2}) and (\ref{eqnestimateTrace}) show
that the first term is bounded by $n^{\o(1)-1}$. Moreover, $\det(\Id
+B)$ is a probability,
so bounded by 1, and the estimates of $\Tr A$ and $\Tr B$, formula
(\ref
{eqnestimateTrace}), show that $e^{\Tr A-\Tr B}-1$ is $\OO(\log n/n)$,
completing the proof.
 \end{pf}
%
%
\begin{lem}\label{leminTheTailGUE}
Let $s(I)=\inf_I\sqrt{4-x^2}$, $\epsilon>0$ and $\alpha_n=\sqrt
{32\log
n}/n$. There are some constants $c_1,c_2>0$ depending only on $\epsilon
$ and $I$,
such that for any $v>\epsilon$ and $n\geq1$,
\[
\E|\{i\dvtx\lambda_i\in I, \lambda_{i+1}-\lambda
_i>(1+v)\alpha
_n/s(I)\}|\leq c_1n^{-c_2 v}.
\]
\end{lem}
\begin{pf}
The first step allows reasoning on fixed intervals instead of gaps.
More precisely, note $I=[a,b]$ with $a<b$, and consider the intervals
of length $(1+\frac{v}{2})\alpha_n/s(I)$ by
successive slips of size $s_n=\frac{v}{2}\alpha_n/s(I)$.
\[
J_k=\biggl[a+k s_n,a+k s_n+\biggl(1+\frac{v}{2}\biggr)\frac{\alpha
_n}{s(I)}\biggr],\qquad 0\leq k\leq p_n=\lfloor(b-a)/s_n\rfloor.
\]
There an injective map associating to any eigenvalues gap of size at
least $(1+v)\alpha_n/s(I)$ an interval $J_k$ included in this gap, for
example, the one with lower index.
Consequently,
\[
\E|\{i\dvtx\lambda_i\in I, \lambda_{i+1}-\lambda
_i>(1+v)\alpha
_n/s(I)\}|\leq
\sum_{k=0}^{p_n}\Prob(J_k=0),
\]
where we use the abbreviation $\Prob(J_k=0)=\Prob
(\lambda_i\notin J_k,1\leq i\leq n)$.
The second step consists in obtaining uniform upper bounds for these
empty intervals probabilities.
For this purpose, the negative correlation property, Lemma~\ref{lemNegCor}, is used by partitioning the interval $J_k$:
\[
J_k=J_k^{(1)}\sqcup J_k^{(2)}\bigsqcup_{j=1}^{q_v} L_k^{(j)}\sqcup M_k,
\]
where:
\begin{itemize}
\item$J_k^{(1)}$ has length $(1-\epsilon')\alpha_n/s(I)$, with
$\epsilon'$ to be chosen positive but sufficiently small compared to
$\epsilon$;
\item$J_k^{(2)}$ has length $(\frac{\epsilon}{2}+\epsilon')\alpha_n/s(I)$;
\item the intervals $L_k^{(j)}$ all have length $\alpha_n/(2s(I))$, and
their number is $q_v=\lfloor v-\epsilon\rfloor$;
\item$M_k$ is a residual interval.
\end{itemize}
The negative correlation property yields
\[
\Prob^{\Gue(n)}(J_k=0)\leq\Prob^{\Gue(n)}
\bigl(J_k^{(1)}=0\bigr)\Prob^{\Gue(n)}\bigl(J_k^{(2)}=0\bigr)\prod
_{j=1}^{q_v}\Prob^{\Gue(n)}\bigl(L_k^{(j)}=0\bigr).
\]
This now can be upper-bounded using Lemma~\ref{lemcompFredholm}
because all intervals are shorter than $\alpha_n/s(I)$: for example,
noting $x$ one extremity of $J_k^{(1)}$,
\begin{eqnarray*}
\Prob^{\Gue(n)}\bigl(J_k^{(1)}=0\bigr)&=&\Prob^{\U(n)}
\bigl(\bigl[0,(1-\epsilon')\alpha_n\sqrt{4-x^2}/s(I)\bigr]=0\bigr)+n^{\o
(1)-1}\\
&\leq&\Prob^{\U(n)}\bigl([0,(1-\epsilon')\alpha_n]=0\bigr)+n^{\o
(1)-1}=n^{-(1-\epsilon')^2+\o(1)}
\end{eqnarray*}
with $\o(1)$ not depending on the index $k$, and where the last
estimate relies on (\ref{eqnest1}). In the same way,
\[
\Prob^{\Gue(n)}\bigl(J_k^{(2)}=0\bigr)\leq n^{-(
{\epsilon
}/{2}+\epsilon')^2+\o(1)}, \qquad\Prob^{\Gue(n)}
\bigl(L_k^{(j)}=0
\bigr)\leq n^{-{1}/{4}+\o(1)}.
\]
Gathering all these results,
\[
\E|\{i\dvtx\lambda_i\in I, \lambda_{i+1}-\lambda
_i>(1+v)\alpha
_n/s(I)\}|\leq n^{1-(1-\epsilon')^2-(
{\epsilon
}/{2}+\epsilon')^2-({1}/{4}+\o(1))q_v}.
\]
We can chose $\epsilon'$ sufficiently small such that $1-(1-\epsilon
')^2-(\frac{\epsilon}{2}+\epsilon')^2<0$ (e.g.,
$\epsilon
'=\epsilon^2/8$). For such a choice, for sufficiently large $n$ not
depending on $v$,
the above exponent is smaller than $-c-\frac{1}{8}\lfloor v-\epsilon
\rfloor$ for some $c>0$. Such a function is smaller than $-c_2 v$ on
$v>\epsilon$, for some $c_2>0$.
\end{pf}
%
%
\begin{lem}\label{lemexpectGUE}
Let $I=[a,b]$ with $a<b$, $s(I)=\inf_I\sqrt{4-x^2}$, $\epsilon\in(0,1)$
and $\alpha_n=\sqrt{32\log n}/n$. Consider a maximal number
\[
\tilde m_n=\bigl\lfloor(b-a)/\bigl((1-\epsilon)\alpha_n/s(I)\bigr)\bigr\rfloor
\]
of disjoint intervals $I_1,\ldots,I_{m_n}$ of length $(1-\epsilon
)\alpha
_n/s(I)$ included in $I$. Let
\[
\tilde M_n=|\{1\leq j\leq m_n \dvtx\lambda_i\notin
I_j,1\leq
i\leq n\}|
\]
be the number of those intervals containing no eigenvalues. Then
\[
\E(\tilde M_n)\gg n^{\delta}
\]
for some $\delta>0$ depending only on $\epsilon$ and $I$.
\end{lem}
\begin{pf}
We first take a restricted interval to avoid the fluctuations in the
spectral measure:
there is $I'=(a',b')\subset I$ with $a'<b'$ such that
%
%
\begin{equation}\label{eqnrestrictedI}
(1-\epsilon)\frac{\sup_{I'}\sqrt{4-x^2}}{s(I)}\leq1-\frac
{\epsilon}{2}.
\end{equation}
Then, for a given interval $I_k=[x,x+(1-\epsilon)\alpha_n/s(I)]$
included in $I'$, from Lem\-ma~\ref{lemcompFredholm},
\begin{eqnarray*}
&&\Prob^{\Gue(n)}(\lambda_i\notin I_k,1\leq i\leq n
)\\
&&\qquad=\Prob^{\U
(n)}\bigl(\theta_i\notin\bigl[0,(1-\epsilon)\sqrt{4-x^2}/s(I)\bigr],1\leq
i\leq
n\bigr)+n^{\o(1)-1}
\end{eqnarray*}
with $\o(1)$ uniform in $x$. From (\ref{eqnrestrictedI}) this is
greater than
\[
\Prob^{\U(n)}\biggl(\theta_i\notin\biggl[0,\biggl(1-\frac
{\epsilon
}{2}\biggr)\biggr],1\leq i\leq n\biggr)+n^{\o(1)-1}=n^{-
(1-{\epsilon}/{2})^2+\o(1)}
\]
from (\ref{eqnest1}). There are $n^{1+\o(1)}$ intervals $I_k$'s
included in $I'$, so as $n\to\infty$
\[
\E(\tilde M_n)\geq n^{1-(1-{\epsilon}/{2})^2+\o
(1)}\gg
n^\delta
\]
with $\delta=(1-(1-\epsilon/2)^2)/2$, for example.
\end{pf}

\subsection{The negative association property}\label{subsecnegAss}
As previously noted, to deduce the asymptotics of the largest gaps, the
correlation between distinct gaps is required.

In the context of a point process on a finite set $\mathcal{E}$, let
$\Lambda$ and $\Lambda'$
be distinct disjoint subsets of $\mathcal{E}$, and write $0^{\Lambda}$
for the event that
the elements of $\Lambda$ are free of particles.
Shirai and Takahashi~\cite{ST} showed that for determinantal point
processes, the empty sets events are negatively correlated.
%
%
\begin{equation}\label{eqnnegDiscr}
\Prob(0^{\Lambda\cup\Lambda'})\leq\Prob(0^{\Lambda})\Prob
(0^{\Lambda'}).
\end{equation}
This negative association property has received considerable attention
in the past few years, in the context of ASEP, for example.
Still, for discrete determinantal point processes,
formula (\ref{eqnnegDiscr}) was generalized to all increasing events
\cite{Lyons},
and general criteria for the negative association property were given
in~\cite{Liggett}.

The following continuous analog of (\ref{eqnnegDiscr}) holds. It can
be shown by a simple discretization, relying on results from
\cite{Liggett,Lyons,ST}. We give another justification,
which relies on a work of Georgii and Yoo~\cite{GY}.
%
%
\begin{lem}[(Negative correlation of the vacuum events)]\label{lemNegCor}
Let $\xi^{(n)}$ be the point process associated to the eigenvalues of
Haar distributed unitary matrix (resp., an element of the $\Gue$). Let
$I_1$ and $I_2$ be compact disjoint subsets of $[0,2\pi)$
(resp.,~$\RR$). Then
%
%
\begin{equation}\label{eqnnegCorr}
\Prob\bigl(\xi^{(n)}(I_1\cup I_2)=0\bigr)\leq
\Prob\bigl(\xi^{(n)}(I_1)=0\bigr)
\Prob\bigl(\xi^{(n)}(I_2)=0\bigr).
\end{equation}
\end{lem}
\begin{pf}
The general negative correlation result, Corollary 3.3 in~\cite{GY},
requires a locally trace class operator $K$
with a restriction on its spectrum. More precisely, consider the
unitary case, the $\Gue$ proof, being similar.
For an operator $K$ acting on $L^2((0,2\pi),\mu)$, if $\spec
(K)\subset[0,1]$,
there exist a unique determinantal point process $\xi$ with kernel $K$;
see~\cite{SoshSurvey}; under the additional hypothesis
$\spec(K)\subset[0,1)$, Georgii and Yoo proved that for any
compact and disjoint Borel sets
$\Lambda\subset\Delta\subset[0,2\pi)$,
\[
\Prob_\mu\bigl(\xi(\Lambda)=0|\xi(\Delta\slash\Lambda
)=0\bigr)\leq
\Prob_\mu\bigl(\xi(\Lambda)=0\bigr).
\]
In the case of Haar distributed unitary matrices, $K=K^{\U(n)}$ is a
nuclear operator with kernel
$K^{\U(n)}(x,y)=\frac{1}{2\pi}\frac{\sin(n(x-y)/2)}{\sin((x-y)/2)}$,
and defines a projection: $1$ is in its spectrum, and the general
statement does not directly apply. To care for this minor problem, look
at the
restriction $K^{\U(n)}_\Lambda$ of $K^{\U(n)}$ to a compact subset
$\Lambda$ of $(0,2\pi)$ ($K^{\U(n)}_\Lambda=P_\Lambda K^{\U(n)}
P_\Lambda$, $P_\Lambda$ being the projection on $\Lambda$). Suppose
that the set $(0,2\pi)\slash\Lambda$ has a nonempty interior.
As a projection of $K^{\U(n)}$, $K^{\U(n)}_\Lambda$ is still
nonnegative and trace class.
As for any determinantal point process,
%
%
\begin{equation}\label{eqnfredDet}
\E\bigl(z^{\xi^{(n)}(\Lambda)}\bigr)=\det\bigl(\Id+(z-1)K^{\U
(n)}_{\Lambda}\bigr)
\end{equation}
in the sense of Fredholm determinants of a trace class operator.
Suppose that $K^{\U(n)}_{\Lambda}$ has an eigenvalue
$\lambda\geq1$.
Then by choosing $z=1-1/\lambda\geq0$, (\ref{eqnfredDet}) yields
\[
\E\bigl(z^{\xi^{(n)}(I)}\bigr)=0\qquad\mbox{if }\lambda>1,\qquad
\Prob\bigl(\xi^{(n)}(I)=0\bigr)=0\qquad\mbox{if }\lambda=1.
\]
In each case, this is absurd because the joint law of the eigenvalues
is absolutely continuous with respect to the Lebesgue measure on
$(0,2\pi)^n$, and $(0,2\pi)\slash\Lambda$ has a nonempty interior:
both quantities need to be strictly positive. Hence $K^{\U
(n)}_{I_1\cup
I_2}$ is a trace class operator with spectrum in $[0,1)$, and the
result from
\cite{GY} applies.
\end{pf}
\begin{rem*}
Another way to prove Lemma~\ref{lemNegCor} is as follows: the
inequality can be stated for the determinantal point process
with kernel $\alpha K^{\U(n)}$ with $0<\alpha<1$, and then the
inequality remains true for $K^{\U(n)}$ by continuity of the application
$K\mapsto\det(\Id-K)$ in the set of trace class operators, by (\ref
{eqnineqTrace}).
\end{rem*}

\section*{Acknowledgments}

This work was initiated while P. Bourgade was visiting the Courant Institute,
and continued in T\'el\'ecom ParisTech: he wishes to thank both
institutions for their kind hospitality and support. Both authors want
to thank Percy Deift for his help and very useful discussions.


%

\printaddresses


\begin{thebibliography}{37}

\bibitem{AGZ}
%
\begin{bbook}[mr]
\bauthor{\bsnm{Anderson},~\bfnm{Greg~W.}\binits{G.~W.}},
\bauthor{\bsnm{Guionnet},~\bfnm{Alice}\binits{A.}} \AND
\bauthor{\bsnm{Zeitouni},~\bfnm{Ofer}\binits{O.}}
(\byear{2010}).
\btitle{An Introduction to Random Matrices}.
\bseries{Cambridge Studies in Advanced Mathematics}
\bvolume{118}.
\bpublisher{Cambridge Univ. Press}, \baddress{Cambridge}.
\bid{mr={2760897}}
\bptnote{check year}%
\bptok{imsref}%
\end{bbook}
%
\endbibitem

\bibitem{BG}
%
\begin{bmisc}[auto:STB|2012/02/29|12:31:17]
\bauthor{\bsnm{Ben~Arous},~\bfnm{G.}\binits{G.}} \AND
\bauthor{\bsnm{Guionnet},~\bfnm{A.}\binits{A.}}
(\byear{2011}).
\bhowpublished{Wigner matrices, Chapter 21. In
\textit{The Oxford Handbook on Random Matrix Theory}
(G. Akemann, J. Baik and P. Di Francesco, eds.) 433--451.
Oxford Univ. Press, Oxford.}
\bptok{imsref}%
\end{bmisc}
%
\endbibitem


\bibitem{BK}
%
\begin{bincollection}[mr]
\bauthor{\bsnm{Ben~Arous},~\bfnm{G{\'e}rard}\binits{G.}} \AND
\bauthor{\bsnm{Kuptsov},~\bfnm{Alexey}\binits{A.}}
(\byear{2009}).
\btitle{R{EM} universality for random {H}amiltonians}.
In \bbooktitle{Spin Glasses: Statics and Dynamics}.
\bseries{Progress in Probability}
\bvolume{62}
\bpages{45--84}.
\bpublisher{Birkh\"auser}, \baddress{Basel}.
\bid{mr={2761980}}
\bptnote{check year}%
\bptok{imsref}%
\end{bincollection}
%
\endbibitem

\bibitem{Liggett}
%
\begin{barticle}[mr]
\bauthor{\bsnm{Borcea},~\bfnm{Julius}\binits{J.}},
\bauthor{\bsnm{Br{\"a}nd{\'e}n},~\bfnm{Petter}\binits{P.}} \AND
\bauthor{\bsnm{Liggett},~\bfnm{Thomas~M.}\binits{T.~M.}}
(\byear{2009}).
\btitle{Negative dependence and the geometry of polynomials}.
\bjournal{J. Amer. Math. Soc.}
\bvolume{22}
\bpages{521--567}.
\bid{doi={10.1090/S0894-0347-08-00618-8}, issn={0894-0347}, mr={2476782}}
\bptok{imsref}%
\end{barticle}
%
\endbibitem

\bibitem{Bui}
%
\begin{barticle}[mr]
\bauthor{\bsnm{Bui},~\bfnm{H.~M.}\binits{H.~M.}}
(\byear{2011}).
\btitle{Large gaps between consecutive zeros of the {R}iemann zeta-function}.
\bjournal{J.~Number Theory}
\bvolume{131}
\bpages{67--95}.
\bid{doi={10.1016/j.jnt.2010.07.010}, issn={0022-314X}, mr={2729210}}
\bptnote{check year}%
\bptok{imsref}%
\end{barticle}
%
\endbibitem


\bibitem{DeiftOpenPbs}
%
\begin{bincollection}[mr]
\bauthor{\bsnm{Deift},~\bfnm{Percy}\binits{P.}}
(\byear{2008}).
\btitle{Some open problems in random matrix theory and the theory of integrable
systems}.
In \bbooktitle{Integrable Systems and Random Matrices}.
\bseries{Contemporary Mathematics}
\bvolume{458}
\bpages{419--430}.
\bpublisher{Amer. Math. Soc.}, \baddress{Providence, RI}.
\bid{mr={2411922}}
\bptok{imsref}%
\end{bincollection}
%
\endbibitem

\bibitem{Deift}
%
\begin{barticle}[mr]
\bauthor{\bsnm{Deift},~\bfnm{P.}\binits{P.}},
\bauthor{\bsnm{Its},~\bfnm{A.}\binits{A.}},
\bauthor{\bsnm{Krasovsky},~\bfnm{I.}\binits{I.}} \AND
\bauthor{\bsnm{Zhou},~\bfnm{X.}\binits{X.}}
(\byear{2007}).
\btitle{The {W}idom--{D}yson constant for the gap probability in random matrix
theory}.
\bjournal{J. Comput. Appl. Math.}
\bvolume{202}
\bpages{26--47}.
\bid{doi={10.1016/j.cam.2005.12.040}, issn={0377-0427}, mr={2301810}}
\bptok{imsref}%
\end{barticle}
%
\endbibitem

\bibitem{DKMVZ}
%
\begin{barticle}[mr]
\bauthor{\bsnm{Deift},~\bfnm{P.}\binits{P.}},
\bauthor{\bsnm{Kriecherbauer},~\bfnm{T.}\binits{T.}},
\bauthor{\bsnm{McLaughlin},~\bfnm{K.~T.~R.}\binits{K.~T.~R.}},
\bauthor{\bsnm{Venakides},~\bfnm{S.}\binits{S.}} \AND
\bauthor{\bsnm{Zhou},~\bfnm{X.}\binits{X.}}
(\byear{1999}).
\btitle{Uniform asymptotics for polynomials orthogonal with respect to varying
exponential weights and applications to universality questions in random
matrix theory}.
\bjournal{Comm. Pure Appl. Math.}
\bvolume{52}
\bpages{1335--1425}.
\bid
{doi={10.1002/(SICI)1097-0312(199911)52:11\&lt;1335::AID-CPA1\&gt;3.0.CO;2-1
}, issn={0010-3640}, mr={1702716}}
\bptok{imsref}%
\end{barticle}
%
\endbibitem

\bibitem{DNT}
%
\begin{barticle}[mr]
\bauthor{\bsnm{Deift},~\bfnm{P.}\binits{P.}},
\bauthor{\bsnm{Nanda},~\bfnm{T.}\binits{T.}} \AND
\bauthor{\bsnm{Tomei},~\bfnm{C.}\binits{C.}}
(\byear{1983}).
\btitle{Ordinary differential equations and the symmetric eigenvalue problem}.
\bjournal{SIAM J. Numer. Anal.}
\bvolume{20}
\bpages{1--22}.
\bid{doi={10.1137/0720001}, issn={0036-1429}, mr={0687364}}
\bptok{imsref}%
\end{barticle}
%
\endbibitem

\bibitem{DelYao}
%
\begin{barticle}[mr]
\bauthor{\bsnm{Delyon},~\bfnm{B.}\binits{B.}} \AND
\bauthor{\bsnm{Yao},~\bfnm{J.}\binits{J.}}
(\byear{2006}).
\btitle{On the spectral distribution of {G}aussian random matrices}.
\bjournal{Acta Math. Appl. Sin. Engl. Ser.}
\bvolume{22}
\bpages{297--312}.
\bid{doi={10.1007/s10255-006-0306-7}, issn={0168-9673}, mr={2216482}}
\bptok{imsref}%
\end{barticle}
%
\endbibitem

\bibitem{Diaconis}
%
\begin{barticle}[mr]
\bauthor{\bsnm{Diaconis},~\bfnm{Persi}\binits{P.}}
(\byear{2003}).
\btitle{Patterns in eigenvalues: The 70th {J}osiah {W}illard {G}ibbs lecture}.
\bjournal{Bull. Amer. Math. Soc. (N.S.)}
\bvolume{40}
\bpages{155--178}.
\bid{doi={10.1090/S0273-0979-03-00975-3}, issn={0273-0979}, mr={1962294}}
\bptok{imsref}%
\end{barticle}
%
\endbibitem

\bibitem{Dyson1}
%
\begin{barticle}[mr]
\bauthor{\bsnm{Dyson},~\bfnm{Freeman~J.}\binits{F.~J.}}
(\byear{1962}).
\btitle{Statistical theory of the energy levels of complex systems. {I}}.
\bjournal{J.~Math. Phys.}
\bvolume{3}
\bpages{140--156}.
\bid{issn={0022-2488}, mr={0143556}}
\bptok{imsref}%
\end{barticle}
%
\endbibitem

\bibitem{Ehr}
%
\begin{barticle}[mr]
\bauthor{\bsnm{Ehrhardt},~\bfnm{Torsten}\binits{T.}}
(\byear{2006}).
\btitle{Dyson's constant in the asymptotics of the {F}redholm
determinant of
the sine kernel}.
\bjournal{Comm. Math. Phys.}
\bvolume{262}
\bpages{317--341}.
\bid{doi={10.1007/s00220-005-1493-4}, issn={0010-3616}, mr={2200263}}
\bptok{imsref}%
\end{barticle}
%
\endbibitem

\bibitem{ESY}
%
\begin{barticle}[mr]
\bauthor{\bsnm{Erd{\H{o}}s},~\bfnm{L{\'a}szl{\'o}}\binits{L.}},
\bauthor{\bsnm{Schlein},~\bfnm{Benjamin}\binits{B.}} \AND
\bauthor{\bsnm{Yau},~\bfnm{Horng-Tzer}\binits{H.-T.}}
(\byear{2011}).
\btitle{Universality of random matrices and local relaxation flow}.
\bjournal{Invent. Math.}
\bvolume{185}
\bpages{75--119}.
\bid{doi={10.1007/s00222-010-0302-7}, issn={0020-9910}, mr={2810797}}
\bptnote{check year}%
\bptok{imsref}%
\end{barticle}
%
\endbibitem

\bibitem{GY}
%
\begin{barticle}[mr]
\bauthor{\bsnm{Georgii},~\bfnm{Hans-Otto}\binits{H.-O.}} \AND
\bauthor{\bsnm{Yoo},~\bfnm{Hyun~Jae}\binits{H.~J.}}
(\byear{2005}).
\btitle{Conditional intensity and {G}ibbsianness of determinantal point
processes}.
\bjournal{J. Stat. Phys.}
\bvolume{118}
\bpages{55--84}.
\bid{doi={10.1007/s10955-004-8777-5}, issn={0022-4715}, mr={2122549}}
\bptok{imsref}%
\end{barticle}
%
\endbibitem

\bibitem{GGK}
%
\begin{bbook}[mr]
\bauthor{\bsnm{Gohberg},~\bfnm{Israel}\binits{I.}},
\bauthor{\bsnm{Goldberg},~\bfnm{Seymour}\binits{S.}} \AND
\bauthor{\bsnm{Krupnik},~\bfnm{Nahum}\binits{N.}}
(\byear{2000}).
\btitle{Traces and Determinants of Linear Operators}.
\bseries{Operator Theory: Advances and Applications}
\bvolume{116}.
\bpublisher{Birkh\"auser}, \baddress{Basel}.
\bid{mr={1744872}}
\bptok{imsref}%
\end{bbook}
%
\endbibitem

\bibitem{Gourdon}
%
\begin{bmisc}[auto:STB|2012/02/29|12:31:17]
\bauthor{\bsnm{Gourdon},~\bfnm{X.}\binits{X.}}
(\byear{2004}).
\bhowpublished{Computation of zeros of the Zeta function. Unpublished
manuscript.}
\bptok{imsref}%
\end{bmisc}
%
\endbibitem

\bibitem{HornJohnson}
%
\begin{bbook}[mr]
\bauthor{\bsnm{Horn},~\bfnm{Roger~A.}\binits{R.~A.}} \AND
\bauthor{\bsnm{Johnson},~\bfnm{Charles~R.}\binits{C.~R.}}
(\byear{1985}).
\btitle{Matrix Analysis}.
\bpublisher{Cambridge Univ. Press}, \baddress{Cambridge}.
\bid{mr={0832183}}
\bptok{imsref}%
\end{bbook}
%
\endbibitem

\bibitem{JohSurvey}
%
\begin{bmisc}[auto:STB|2012/02/29|12:31:17]
\bauthor{\bsnm{Johansson},~\bfnm{K.}\binits{K.}}
(\byear{2005}).
\bhowpublished{Random matrices and determinantal processes.
Lecture notes, Summer School on Mathematical Statistical Mechanics
at Ecole de Physique, Les Houches.}
\bptok{imsref}%
\end{bmisc}
%
\endbibitem

\bibitem{Kallenberg}
%
\begin{bbook}[mr]
\bauthor{\bsnm{Kallenberg},~\bfnm{Olav}\binits{O.}}
(\byear{1983}).
\btitle{Random Measures}, \bedition{3rd} ed.
\bpublisher{Akademie Verlag}, \baddress{Berlin}.
\bid{mr={0818219}}
\bptok{imsref}%
\end{bbook}
%
\endbibitem

\bibitem{KatzSarnak1}
%
\begin{bbook}[mr]
\bauthor{\bsnm{Katz},~\bfnm{Nicholas~M.}\binits{N.~M.}} \AND
\bauthor{\bsnm{Sarnak},~\bfnm{Peter}\binits{P.}}
(\byear{1999}).
\btitle{Random Matrices, {F}robenius Eigenvalues, and Monodromy}.
\bseries{American Mathematical Society Colloquium Publications}
\bvolume{45}.
\bpublisher{Amer. Math. Soc.}, \baddress{Providence, RI}.
\bid{mr={1659828}}
\bptok{imsref}%
\end{bbook}
%
\endbibitem

\bibitem{KeaSna}
%
\begin{barticle}[mr]
\bauthor{\bsnm{Keating},~\bfnm{J.~P.}\binits{J.~P.}} \AND
\bauthor{\bsnm{Snaith},~\bfnm{N.~C.}\binits{N.~C.}}
(\byear{2000}).
\btitle{Random matrix theory and {$\zeta(1/2+it)$}}.
\bjournal{Comm. Math. Phys.}
\bvolume{214}
\bpages{57--89}.
\bid{doi={10.1007/s002200000261}, issn={0010-3616}, mr={1794265}}
\bptok{imsref}%
\end{barticle}
%
\endbibitem

\bibitem{Kra}
%
\begin{barticle}[mr]
\bauthor{\bsnm{Krasovsky},~\bfnm{I.~V.}\binits{I.~V.}}
(\byear{2004}).
\btitle{Gap probability in the spectrum of random matrices and
asymptotics of
polynomials orthogonal on an arc of the unit circle}.
\bjournal{Int. Math. Res. Not. IMRN}
\bvolume{25}
\bpages{1249--1272}.
\bid{doi={10.1155/S1073792804140221}, issn={1073-7928}, mr={2047176}}
\bptok{imsref}%
\end{barticle}
%
\endbibitem

\bibitem{Lyons}
%
\begin{barticle}[mr]
\bauthor{\bsnm{Lyons},~\bfnm{Russell}\binits{R.}}
(\byear{2003}).
\btitle{Determinantal probability measures}.
\bjournal{Publ. Math. Inst. Hautes \'Etudes Sci.}
\bvolume{98}
\bpages{167--212}.
\bid{doi={10.1007/s10240-003-0016-0}, issn={0073-8301}, mr={2031202}}
\bptok{imsref}%
\end{barticle}
%
\endbibitem

\bibitem{Mehta}
%
\begin{bbook}[mr]
\bauthor{\bsnm{Mehta},~\bfnm{Madan~Lal}\binits{M.~L.}}
(\byear{2004}).
\btitle{Random Matrices},
\bedition{3rd} ed.
\bseries{Pure and Applied Mathematics (Amsterdam)}
\bvolume{142}.
\bpublisher{Elsevier}, \baddress{Amsterdam}.
\bid{mr={2129906}}
\bptok{imsref}%
\end{bbook}
%
\endbibitem

\bibitem{Montg}
%
\begin{bincollection}[mr]
\bauthor{\bsnm{Montgomery},~\bfnm{H.~L.}\binits{H.~L.}}
(\byear{1973}).
\btitle{The pair correlation of zeros of the zeta function}.
In \bbooktitle{Analytic Number Theory ({P}roc. {S}ympos. {P}ure
{M}ath., {V}ol.
{XXIV}, {S}t. {L}ouis {U}niv., {S}t. {L}ouis, {M}o., 1972)}
\bpages{181--193}.
\bpublisher{Amer. Math. Soc.}, \baddress{Providence, RI}.
\bid{mr={0337821}}
\bptok{imsref}%
\end{bincollection}
%
\endbibitem

\bibitem{Odl}
%
\begin{barticle}[mr]
\bauthor{\bsnm{Odlyzko},~\bfnm{A.~M.}\binits{A.~M.}}
(\byear{1987}).
\btitle{On the distribution of spacings between zeros of the zeta function}.
\bjournal{Math. Comp.}
\bvolume{48}
\bpages{273--308}.
\bid{doi={10.2307/2007890}, issn={0025-5718}, mr={0866115}}
\bptok{imsref}%
\end{barticle}
%
\endbibitem

\bibitem{ST}
%
\begin{barticle}[mr]
\bauthor{\bsnm{Shirai},~\bfnm{Tomoyuki}\binits{T.}} \AND
\bauthor{\bsnm{Takahashi},~\bfnm{Yoichiro}\binits{Y.}}
(\byear{2003}).
\btitle{Random point fields associated with certain {F}redholm determinants.
{II}. {F}ermion shifts and their ergodic and {G}ibbs properties}.
\bjournal{Ann. Probab.}
\bvolume{31}
\bpages{1533--1564}.
\bid{doi={10.1214/aop/1055425789}, issn={0091-1798}, mr={1989442}}
\bptok{imsref}%
\end{barticle}
%
\endbibitem

\bibitem{Sosh1}
%
\begin{barticle}[mr]
\bauthor{\bsnm{Soshnikov},~\bfnm{Alexander}\binits{A.}}
(\byear{1998}).
\btitle{Level spacings distribution for large random matrices: {G}aussian
fluctuations}.
\bjournal{Ann. of Math. (2)}
\bvolume{148}
\bpages{573--617}.
\bid{doi={10.2307/121004}, issn={0003-486X}, mr={1668559}}
\bptok{imsref}%
\end{barticle}
%
\endbibitem

\bibitem{SoshSurvey}
%
\begin{barticle}[mr]
\bauthor{\bsnm{Soshnikov},~\bfnm{A.}\binits{A.}}
(\byear{2000}).
\btitle{Determinantal random point fields}.
\bjournal{Russian Math. Surveys}
\bvolume{55}
\bpages{923--975}.
\bptok{imsref}%
\end{barticle}
%
\endbibitem

\bibitem{Sosh2}
%
\begin{barticle}[mr]
\bauthor{\bsnm{Soshnikov},~\bfnm{Alexander}\binits{A.}}
(\byear{2005}).
\btitle{Statistics of extreme spacing in determinantal random point processes}.
\bjournal{Mosc. Math. J.}
\bvolume{5}
\bpages{705--719, 744}.
\bid{issn={1609-3321}, mr={2241818}}
\bptok{imsref}%
\end{barticle}
%
\endbibitem

\bibitem{Szego}
%
\begin{bbook}[mr]
\bauthor{\bsnm{Szeg{\H{o}}},~\bfnm{G{\'a}bor}\binits{G.}}
(\byear{1975}).
\btitle{Orthogonal Polynomials}, \bedition{4th} ed.
\bseries{American Mathematical Society Colloquium Publications}
\bvolume{XXIII}.
\bpublisher{Amer. Math. Soc.}, \baddress{Providence, RI}.
\bid{mr={0372517}}
\bptok{imsref}%
\end{bbook}
%
\endbibitem

\bibitem{TV}
%
\begin{barticle}[mr]
\bauthor{\bsnm{Tao},~\bfnm{Terence}\binits{T.}} \AND
\bauthor{\bsnm{Vu},~\bfnm{Van}\binits{V.}}
(\byear{2011}).
\btitle{Random matrices: Universality of local eigenvalue statistics}.
\bjournal{Acta Math.}
\bvolume{206}
\bpages{127--204}.
\bid{doi={10.1007/s11511-011-0061-3}, issn={0001-5962}, mr={2784665}}
\bptnote{check year}%
\bptok{imsref}%
\end{barticle}
%
\endbibitem

\bibitem{Vinson}
%
\begin{bmisc}[auto:STB|2012/02/29|12:31:17]
\bauthor{\bsnm{Vinson},~\bfnm{J.}\binits{J.}}
(\byear{2001}).
\bhowpublished{Closest spacing of eigenvalues. Ph.D.
thesis, Princeton Univ.}
\bptok{imsref}%
\end{bmisc}
%
\endbibitem

\bibitem{Widom}
%
\begin{barticle}[mr]
\bauthor{\bsnm{Widom},~\bfnm{Harold}\binits{H.}}
(\byear{1971}).
\btitle{The strong {S}zeg{\H o} limit theorem for circular arcs}.
\bjournal{Indiana Univ. Math. J.}
\bvolume{21}
\bpages{277--283}.
\bid{issn={0022-2518}, mr={0288495}}
\bptok{imsref}%
\end{barticle}
%
\endbibitem

\end{thebibliography}
\end{document}